\newif\ifdraft
\draftfalse

\ifdraft
\documentclass[journal,12pt,onecolumn,draftclsnofoot]{IEEEtran}
\usepackage[mathlines]{lineno}
\modulolinenumbers[5]
\linenumbers\relax 
\else
\documentclass[journal]{IEEEtran}
\fi

\usepackage{amsmath}                 
\usepackage{amssymb}                 
\usepackage{mathrsfs}                
\usepackage{graphicx}

\usepackage{algorithm}               
\usepackage{algpseudocode}           
\algrenewcommand\algorithmicrequire{\textbf{Input:}}
\algrenewcommand\algorithmicensure {\textbf{Output:}}

\usepackage[table,xcdraw]{xcolor}
\usepackage{booktabs}               
\usepackage{etoolbox}
\makeatletter \patchcmd{\@makecaption} {\scshape} {} {} {} \makeatother
\usepackage{color}                  

\usepackage[square,numbers,sort&compress]{natbib} 
\usepackage[capitalise]{cleveref}                 
\crefname{equation}{}{}                           
\Crefname{equation}{}{}                           

\newcommand{\cost}{\cos\theta}
\newcommand{\db}{\,\mbox{dB}}

\newcommand{\Fpp}{\mathcal{F}_\text{pp}}
\newcommand{\F}{\mathcal{F}}

\newcommand{\ra}{\rightarrow}

\newcommand{\norm}[1]{\left\Vert{#1}\right\Vert}

\newcommand{\Q}{\mathcal{Q}}

\newcommand{\Spp}{\mathcal{S}^*_\text{pp}}
\newcommand{\Rpp}{\mathcal{R}_\text{pp}}
\newcommand{\A}{\mathcal{A}}
\newcommand{\D}{\mathcal{D}}
\newcommand{\M}{\mathcal{M}}
\newcommand{\N}{\mathcal{N}}
\renewcommand{\O}{\mathcal{O}}
\newcommand{\R}{\mathbb{R}}
\newcommand{\W}{\mathcal{W}}
\renewcommand{\S}{\mathcal{S}^*}
\newcommand{\sint}{\sin\theta}
\renewcommand{\v}[1]{{\bf{#1}}}
\newcommand{\dint}[4]{{\int_{#1}^{#2}#3\,\text{d}#4}}
\newcommand{\SNR}[1]{\text{SNR}_{\text{dB}}\left(#1\right)}

\author{
	Shahar Tsiper and Yonina C. Eldar
	
	\thanks{This project has received funding from the European Union's Horizon 2020 research and innovation program under grant No. 646804-ERC-COG-BNYQ.\@
		
		S. Tsiper (tsiper@technion.ac.il) and Y. C. Eldar (yonina@ee.technion.ac.il) are with the Department of Electrical Engineering,
		Technion --- Israel Institute of Technology, Haifa 32000.}
}

\newtoggle{TMI}
\togglefalse{TMI}

\title{RAPToR:\ A Resampling Algorithm for Pseudo-Polar based Tomographic Reconstruction}

\begin{document}

\maketitle
\begin{abstract}
We propose a stable and fast reconstruction technique for parallel-beam (PB) tomographic X-ray imaging, relying on the discrete pseudo-polar (PP) Radon transform.
Our main contribution is a resampling method, based on modern sampling theory, that transforms the acquired PB measurements to a PP grid.
The resampling process is both fast and accurate, and in addition, simultaneously denoises the measurements, by exploiting geometrical properties of the tomographic scan.
The transformed measurements are then reconstructed using an iterative solver with total variation (TV) regularization.
We show that reconstructing from measurements on the PP grid, leads to improved recovery, due to the inherent stability and accuracy of the PP Radon transform, compared with the PB Radon transform.
We also demonstrate recovery from a reduced number of PB acquisition angles, and high noise levels. Our approach is shown to achieve superior results over other state-of-the-art solutions, that operate directly on the given PB measurements.
The proposed method can benefit fan-beam and/or cone-beam projections by coupling it with a rebinning process.
\end{abstract}

\ifdraft\else
\fi

\section{Introduction}\label{sec:intro}

Projectional radiography imaging is the most widely used imaging modality in medicine today.
The computed tomography (CT) technique for 3D imaging was invented in the early 1970's by Sir Godfrey Newbold Hounsefield, and since then CT scanners have come a long way.
Today's scanners can perform full body scans in under 30 seconds, incorporating detector arrays with hundreds of elements that acquire more than a thousand different readings in every revolution of the scanner.
Despite significant technological progress in the manufacturing of CT scanners, the reconstruction algorithms in use have only slightly evolved from the very first methods developed more than 30 years ago.
The vast majority of commercial algorithms are based on the filtered back projection (FBP) algorithm~\cite{Hsieh2003,Kak1988c} for parallel-beam (PB) and fan-beam scans, or the Feldkamp-Davis-Kress (FDK) algorithm~\cite{Feldkamp1984,Hsieh2003} for cone-beam/helical scans.
Modern reconstruction approaches, which rely on iterative optimization based solvers, are seldom used in commercial machines due to the prohibitive computational complexity involved, leading to long reconstruction durations.


Recently, a paradigm for processing discrete tomographic measurements has been proposed, which relies on the pseudo-polar (PP) Radon transform (PPRT)~\cite{Averbuch2001,Averbuch2008,Averbuch2008a}.
PPRT is an algebraic mapping that relates an object to its discrete tomographic projections.
This transform has many advantages, further described below, that make it an appealing framework for reconstruction of tomographic images.
However, the PP framework assumes the projections are taken over a non-uniform set of angles and detectors, and therefore cannot be used directly on today's tomographic scans.

There are three main advantages to using the PP framework.
First, the PP grid points are closer to the polar grid points than the Cartesian grid, thus interpolating polar measurements to the PP grid is easier and more accurate than transforming them to a Cartesian grid, which is the standard practice today.
Second, the PP Fourier transform (PPFT) and the PPRT can be computed with a fast accurate algorithm~\cite{Averbuch2008,Averbuch2008a}.
The adjoint PPRT operation is computable at the same speed and accuracy as the forward PPFT, thus allowing a fast iterative inversion scheme.
The third advantage, is that the algebraic system that describes the PPRT has a significantly lower condition number than an equivalent PB system.
This allows to reconstruct from the PP measurements significantly faster and more accurately than from PB samples, with lower noise amplification.
Although these are significant advantages, the PP framework is not used today for commercial tomographic reconstruction, since CT machines do not acquire measurements on this grid.

\begin{figure}[t]
    \centering
    \includegraphics[width=\columnwidth]{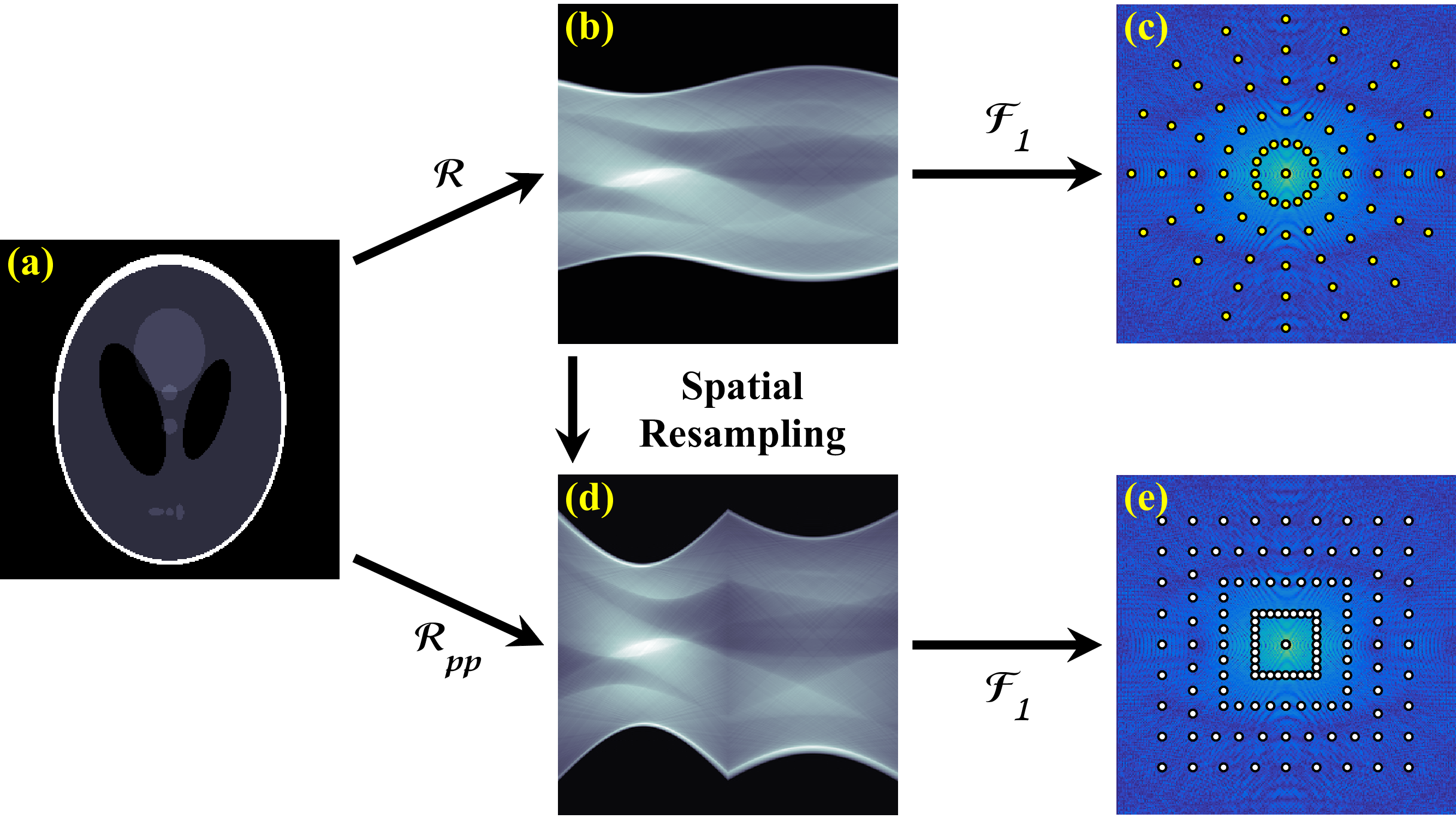}
    \caption{
        \textbf{Diagram of the different transforms.} (a) The Shepp-Logan phantom $f \left( x,y \right) $. (b) The PB sinogram $p(t,\theta)$ defined in \cref{eq:RadonTransform} and \cref{eq:PB_Sampling}. (c) The 1D-DFT of the PB sinogram $\v s$ \cref{eq:PB_1D_DFT}, which samples the 2D-CSFT of $f$ in polar coordinates. (d) The PP sinogram $\hat{\v p}$ \cref{eq:1D_IDFT_p_Sum_s}. (e) The 1D-DFT of the PP sinogram $\hat{\v s}$ \cref{eq:PPFT_Definition}, which evaluates the 2D-CSFT of $f$ on the PP grid.\label{fig:Transforms}}
\end{figure}

In this work, we propose a way to bridge the gap between PP theory and real-world CT scans, using a Resampling Algorithm for Pseudo-polar based Tomographic Reconstruction (RAPToR).
Our main contribution is a resampling algorithm, which is based on modern sampling tools, that transforms actual CT measurements to comply with the PP model.
This approach relies on a mathematical derivation of a subspace, in which the tomographic measurements lie.
The subspace and its respective kernel take into account the geometrical structure of the tomographic scan, expanding the work in~\cite{Rattey1981}, and utilize tools from generalized sampling theory~\cite{Eldar2015}.
The resampling technique presented is shown to be numerically accurate, and to simultaneously denoise the measurements, while being computationally efficient.
We demonstrate that a fast iterative solver (FISTA-TV~\cite{Beck2009b,DanielP.Palomar2010}) coupled with the PP framework, produces excellent reconstructions from the transformed PP measurements, surpassing other state-of-the-art solutions.
Better performance is maintained for both low signal to noise ratios (SNR) and when given a reduced number of measurements.
These translate to radiation dosage reduction for CT scanning.

Prior works have already suggested ways to benefit from the PP framework.
The authors in~\cite{Lee2008,Fahimian2010,Tomography2010} proposed iterative reconstruction schemes based on the PP transform, yet assumed that the measurements were acquired directly on the PP grid.
Unfortunately, this assumption is not valid for CT scanners in use today.
In~\cite{Averbuch2001} the uniform PB measurements are interpolated to the non-uniform PP grid.
The suggested interpolation process is performed in the frequency domain, and involves extrapolations, which suffer from numerical inaccuracy, usually leading to unwanted reconstruction artifacts.
Furthermore, no results or comparisons were provided to assess the quality of the interpolation.
In~\cite{Hashemi2014} the measurements are acquired over uniformly spaced detectors, and non-uniform angles, so that the angles comply with the PP paradigm, which might be challenging for commercial CTs.
The measurements are then transformed to the PP grid, by performing frequency based interpolation solely over the detector axis.

Our approach offers a fundamentally different methodology that generates PP measurements from a PB projection, by exploiting the inherent structural properties of the tomographic scan.
All the interpolations in our solution are performed in the spatial domain of the measurements, and operate with a low computational burden.
In addition, we do not assume that the acquisition angles comply with the PP angles.
The computational complexity of both resampling and reconstruction is on the order of $\O(N^2\log N)$, where $N^2$ is the number of pixels in the reconstructed CT axial scan.
We show that our method outperforms state-of-the-art algorithms for tomographic reconstruction, both in of computational complexity and in quality.

While the method presented in this manuscript focuses on PB measurements, we have found it applicable as well to both fan-beam and cone-beam reconstructions, as acquired in modern CT devices.
Our method can be used on these acquisition schemes by first applying a rebinning step~\cite{Kudo1999,Taguchi2000,Schaller2000,Louis2006} on the fan-beam/cone-beam measurements that transforms them into PB measurements. 
Rebinning to PB measurements is used today in the reconstruction flow of commercial devices for enabling fast filtration and back-projecting using the PB FBP algorithm.
Our approach can then replace the traditional FBP step and enable fast reconstruction via the PP grid, with improved noise resistance and accuracy, even when compared to direct reconstruction from the fan/cone-beam measurements.
This application-oriented scheme for reconstructing from fan and cone-beam acquisitions is further discussed in the 
\iftoggle{TMI}{Supplementary Information}{Appendix~\ref{sec:Fan_Appendix}}.


The paper is organized as follows.
\Cref{sec:Problem} contains the mathematical prerequisites and defines the problem.
In \cref{sec:Pseudo_Polar} we review the PP transform and discuss its properties, followed by \cref{sec:Subspace} where we present our subspace approach for resampling the sinogram.
\Cref{sec:Reconstruction} suggests a reconstruction method, composed of resampling to the PP grid and then reconstructing using a modern and a fast iterative solver, regularized by the total variation (TV)~\cite{Rudin1992} norm.
An experiment is detailed in \cref{sec:Setup}, followed by a conclusion in \cref{sec:Discussion}.

\section{Problem Formulation}\label{sec:Problem}

We begin by describing the acquisition process of a PB tomographic scan, and then formulate the problem we aim to solve, first using continuous notation and then defining the discrete sampling process.
The 2D axial object to be scanned is denoted by $f(x,y):\R^2\ra\R_+$, where $\{x,y\}$ are spatial coordinates.
The value of the scanned object $f$ in each point, represents the X-ray attenuation coefficient.
We assume that the support of the object is confined to a circle with radius $R$, so that $f(x,y)=0$, for $x^2+y^2\geq R^2$.
In addition, we define $\v f\in\R^{N\times N}_+$ as the discrete image on a Cartesian grid, such that $\v f[u,v]=f(uT,vT)$, where $T$ is the spacing between two pixels and $\{u,v\}\in[-\frac{N}{2},\frac{N}{2}-1]$.
Requiring that $2R<NT$ ensures that the discrete image spans the entire support of the object.

The CT measurements are composed from projections, denoted by $p(t,\theta):\R^2\ra\R_+$, where $\theta$ is the projection angle and $t$ is the detector location. These projections are commonly referred to as the sinogram.
The relation between the sinogram and the scanned object is given by the PB Radon transform
\begin{align}\label{eq:RadonTransform}
    p(t,\theta) =&
    \dint{-\infty}{\infty}{\dint{0}{\pi}
    {f(x,y)\delta(t-x\cost-y\sint)}\theta}t\,.
\end{align}
The 1D and 2D continuous space Fourier transforms (CSFT) of the sinogram, and the 2D-CSFT of the object are denoted by $S$, $P$ and $F$ respectively, and given by
\begin{align}\label{eq:1D_Sino_CSFT}
    S\left(\omega_{t},\theta\right)=&
        \dint{-\infty}{\infty}{p\left(t,\theta\right)
        e^{-j\omega_{t}t}}t                   \\
    \label{eq:2D_Sino_CSFT}
    P\left( \omega_t ,\omega_\theta\right)=&
        \dint{-\infty}{\infty}{ S\left( \omega_t,\theta \right) e^{-j\omega_\theta\theta}}{\theta} \\
    \label{eq:2D_Object_CSFT}
    F \left( \omega_x,\omega_y \right) =&
        \dint{-\infty}{\infty}{\dint{-\infty}{\infty}{f\left(x,y\right)
        e^{-j(\omega_x x + \omega_y y)}}x}y\,.
\end{align}
Note that the definition of $P$ in \cref{eq:2D_Sino_CSFT} is not intuitive, since a CSFT performed on an angular variable is not a common procedure.
In \cref{sec:Subspace} we use \cref{eq:2D_Sino_CSFT} for constructing the subspace of the measurements.

The function $S$ and the 2D-CSFT of the object $F$ are related by
\begin{align}\label{eq:FourierSliceThm}
    F(\omega_t\cos\theta,\omega_t\sin\theta) = S(\theta,\omega_t)\,.
\end{align}
This property is commonly referred to as the Fourier slice theorem~\cite{Hsieh2003}.
Its intuitive meaning is that a PB projection of an object taken at an angle $\theta$, is equivalent to a line in the frequency domain, sampling the 2D-CSFT of the object, and intersecting the horizontal frequency axis $\omega_x$ with an angle of $\theta$.

\begin{figure}[t]
    \centering
    \includegraphics[width=\columnwidth]{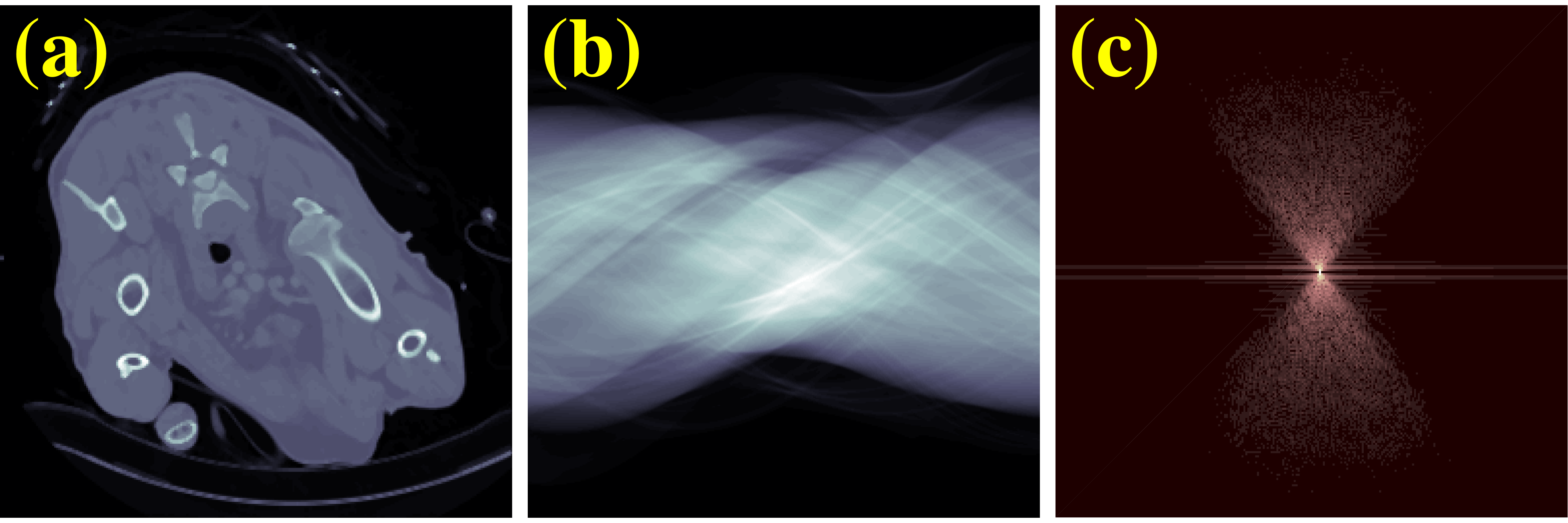}
    \caption{\textbf{Presenting the object and its PB projections.}
        (a) The ground-truth digitized sheep phantom $f(x,y)$, with full dynamic range (b) the polar sinogram $p\left(t,\theta\right)$ \cref{eq:RadonTransform}, (c) the 2D-CSFT of the sinogram $P\left(\omega_t,\omega_\theta\right)$ \cref{eq:2D_Sino_CSFT}.
        The spectral support structure is shown to be limited to the region $\Omega$ \cref{eq:OmegaSS}.
        \label{fig:PolarSinograms}}
\end{figure}

The sampled PB sinogram $\v p\in\R^{N \times \lceil\pi/\Theta\rceil}_+$ is given by
\begin{align} \label{eq:PB_Sampling}
    \v{p}[m,n] & = p(mT,n\Theta) \,.
\end{align}
Here, $\Theta$ and $T$ are the spacings between two adjacent projection angles and detectors respectively.
We use the short-hand notation
\begin{align} \label{eq:S_Definition}
    \v p = \S p
\end{align}
where the PB sampling operator $\S$ operates on a continuous sinogram $p$, producing its discrete samples $\v p$.
The full PB scan takes place over an angular range of $\pi$ radians and a detector axis of length $NT$.
Note that the PB sampling process is uniform with respect to the variables of $\v p$, namely $(t,\theta)$.

The 1D discrete Fourier transform (DFT) of $\v p$ is defined as
\begin{align} \label{eq:PB_1D_DFT}
    \v s[k,n] &= \sum_{m=-N/2}^{N/2-1} \v p[m,n] e^{-2\pi jmk/(N+1) }\,,
\end{align}
and denoted $\v s = \F_1 \v p$.
Here $\F_1$ is the 1D-DFT operator operating on the columns of its argument.
Using \cref{eq:FourierSliceThm,eq:PB_Sampling,eq:PB_1D_DFT}, the values of $F$ can be related to $\v s$ by
\begin{align} \label{eq:PB_Discrete_FST}
    \v s[k,n]= F\left(\frac{k\pi}{NT}\cos(n\Theta) , \frac{k\pi}{NT}\sin(n\Theta)  \right),
\end{align}
for $k\in[-\frac{N}{2},\frac{N}{2}-1]$ and $n\in\left[ 0,\lceil \frac{\pi}{\Theta} \rceil\right]$.
This establishes that a PB scan is equivalent to sampling the 2D-CSFT of the scanned object on a polar grid.
An example of a polar grid in the frequency domain can be seen in \cref{fig:Transforms}(c).

One approach for reconstructing a PB scan is to recover $f$ from the polar measurements of the 2D-CSFT given in \cref{eq:PB_Discrete_FST}, by first transforming the measurements $\v p$ to $\v s$ using \cref{eq:PB_1D_DFT}.
This approach was recently demonstrated by~\cite{Fessler2003,Kiperwas2016}.
However, inverting the 2D-CSFT from a set of points known on a polar grid is a difficult task, both in terms of computational complexity and in terms of inversion accuracy.
We will compare our method to the SParse Uniform ReSampling (SPURS)~\cite{Kiperwas2016} algorithm,
a state-of-the-art algorithm for performing non-uniform to uniform resampling in the frequency domain, that surpassed the non-uniform fast Fourier transform algorithm proposed in~\cite{Fessler2003}.

The main algorithms in use today, FBP~\cite{Kak1988c} and FDK~\cite{Feldkamp1984}, approximate the solution to this inverse problem, yet lead to substantial reconstruction artifacts.
They offer a direct one-shot reconstruction approach, without applying any interpolations, by performing one-dimensional filtering on the sinogram and then back-projecting it.
However, these direct methods require a multitude of measurements to produce an image of clinical quality and produce more artifacts than recent iterative approaches as further discussed in~\cite{Pan2009}. 
In addition, their computational complexity is on the order of $\O(N^3\log N)$, which is higher than that of a single iterative step of our algorithm.

We propose a new reconstruction scheme composed of two steps: we first transform the given PB measurements $\v p$ \cref{eq:PB_Sampling} to a modified PP sinogram $\hat{\v p}$, defined below and illustrated in \cref{fig:Transforms}(d), by the ``spatial resampling'' arrow.
Next, we solve the linear PP system given by $\hat{\v p} = \Rpp\v f$, obtaining the reconstructed object $\hat{\v f}$, where $\Rpp$ is the PPRT linear operator that takes an image  $\v f \in \R_+^{N\times N}$, and returns the PP projection samples $\hat{\v p}$.
By resampling our measurements from the polar grid to the PP grid, shown in \cref{fig:Transforms}(e), we gain all the advantages of the PP framework, including improved numerical accuracy, algebraic stability and fast iterative inversion.
The fast inversion is due to the low condition number associated with the PPRT operator, compared to a PB Radon operator. These advantages eventually lead to more accurate reconstructions while reducing the overall computational complexity.

\section{The Pseudo-Polar Transform}\label{sec:Pseudo_Polar}
In this section we define the PP transform and mathematically relate it to PB scanning.
The relations established here, are used in the next sections for formulating a robust resampling scheme that transforms the measurements between their acquisition grid and the PP grid, alleviating the main drawback of the PP framework.

The PPFT evaluates the 2D discrete space Fourier transform (DSFT) of an $N \times N$ image on a PP grid, seen in \cref{fig:Transforms}{(e)}, over the frequency domain $(\omega_x,\omega_y)$.
This grid is also known as a ``concentric squares'' grid~\cite{Averbuch2001}, and occasionally referred to as ``equally sloped tomography''.
The PPFT operator $\Fpp$ computes the 2D-DTFT values of $\v f$ on the PP grid, and is given by
\begin{align}\label{eq:PPFT_Definition}
    \hat{\v s}[k,n] = \Fpp \v f = \begin{cases}
    \sum\limits_{u,v} \v f[u,v] e^{-\frac{2\pi j}{2N+1} \left(-\frac{2nk}{N}u+kv\right)}, &
        n\in\N_1, \\
    \sum\limits_{u,v} \v f[u,v] e^{-\frac{2\pi j}{2N+1} \left(+ku-\frac{2nk}{N}v\right)}, &
        n\in\N_2,
    \end{cases}
\end{align}
where $\N_1 \triangleq [-\frac{N}{2},\frac{N}{2}-1]$, and $\N_2 \triangleq [\frac{N}{2},\frac{3N}{2}]  $.
We define the PPRT using its connection to the PPFT, by applying an inverse 1D-DFT on the columns of $\hat{\v s}$:
\begin{align}\label{eq:1D_IDFT_p_Sum_s}
    \hat{\v p} &= \Rpp \v f = \F_1^{-1} \hat{\v s} = \F_1^{-1}\Fpp \v f\,,\\
    \hat{\v p}[m,n] &=\sum\limits_{k=-N}^{N} \hat{\v s}[k,n]e^{ \frac{2\pi j}{2N+1}km },
\end{align}
where $m\in[-N,N]$ and $n\in[-\frac{N}{2},\frac{3N}{2}]$, such that ${\hat{\v p}\in \R^{(2N+1)\times(2N+1)}}$ is the PP sampled sinogram.
The relations between the PPFT, PPRT and the PB tomographic scan are illustrated in \cref{fig:Transforms}.

Next, our goal is to obtain a relation between the PP sampled sinogram $\hat{\v p}$ and the continuous sinogram $p$, in order to define the PP sampling operator, which is used for resampling the samples onto the PP grid.
By combining the discrete Fourier slice theorem (Theorem 1 in~\cite{Averbuch2001}), together with \cref{eq:PB_1D_DFT,eq:PB_Discrete_FST,eq:PPFT_Definition,eq:1D_IDFT_p_Sum_s}, we get
\begin{align}\label{eq:PP_Spatial_Sampling}
    \hat{\v p}[m,n] = & \Spp p = p\left(m \hat{T}_n, n\hat{\Theta}_n\right),
\end{align}
where $\Spp$ is the PP sampling operator, and
\begin{align}\label{eq:PP_Tn_Theta_n}
    \left( \hat{T}_n , \hat{\Theta}_n \right) = & \begin{cases}
        \left( \frac{T}{\sqrt{1+{\left(\frac{2n}{N}\right)}^2}} ,\ \tan^{-1}(\frac{2n}{N}) \right),
        & n\in\N_1, \\
        \left( \frac{T}{\sqrt{1+{\left(\frac{2n-N}{N}\right)}^2}} ,\ \cot^{-1}(\frac{2n-N}{N}) \right),
        & n\in\N_2\,    ,
    \end{cases}
\end{align}
define the non-uniform PP grid points in the $(t,\theta)$ domain.
The above expressions \cref{eq:PP_Spatial_Sampling,eq:PP_Tn_Theta_n} define the direct sampling method for the PP sinogram, as can be seen in \cref{fig:Transforms}(d).

To summarize, by relying on \cref{eq:PB_Sampling,eq:PP_Spatial_Sampling,eq:PP_Tn_Theta_n} we formulate in the next section our first step; a subspace-based resampling from the given PB sinogram $\v p$ to a corresponding PP counterpart $\hat{\v p}$.
Once we have the PP measurements, we reconstruct the object $\v f$ using the PP operators given in \cref{eq:PPFT_Definition,eq:1D_IDFT_p_Sum_s}.

\begin{figure}[t]
    \centering
    \includegraphics[width=\columnwidth]{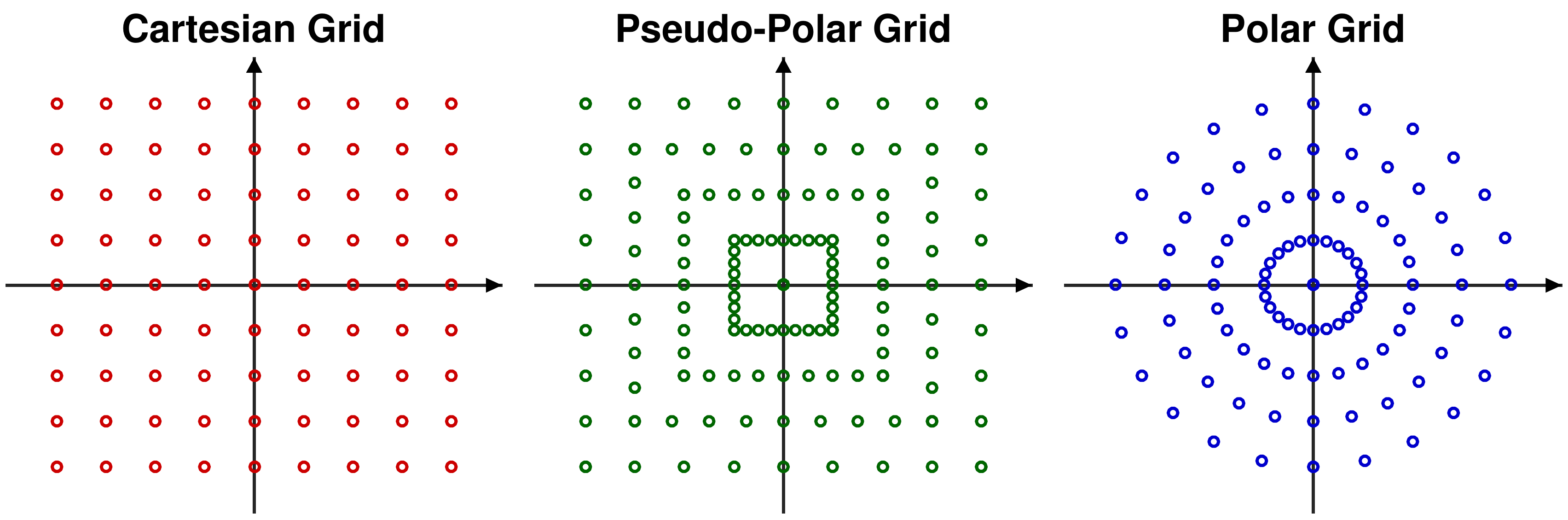}
    \caption{The different acquisition grids in Fourier (frequency) domain.
            The Cartesian grid (left) is the grid used by the 2D-DFT\@.
            After parallel-beam scanning we obtain the values of the Fourier transform of the scanned object on the polar grid (right).
            The pseudo-polar grid (middle), acts as an intermediate step between the two grids.}
    \label{fig:Grids}
\end{figure}

\section{The Sinogram Subspace}
\label{sec:Subspace}

\subsection{Mathematical Definition}

\begin{figure}[t]
    \centering
    \includegraphics[width=0.85\columnwidth]{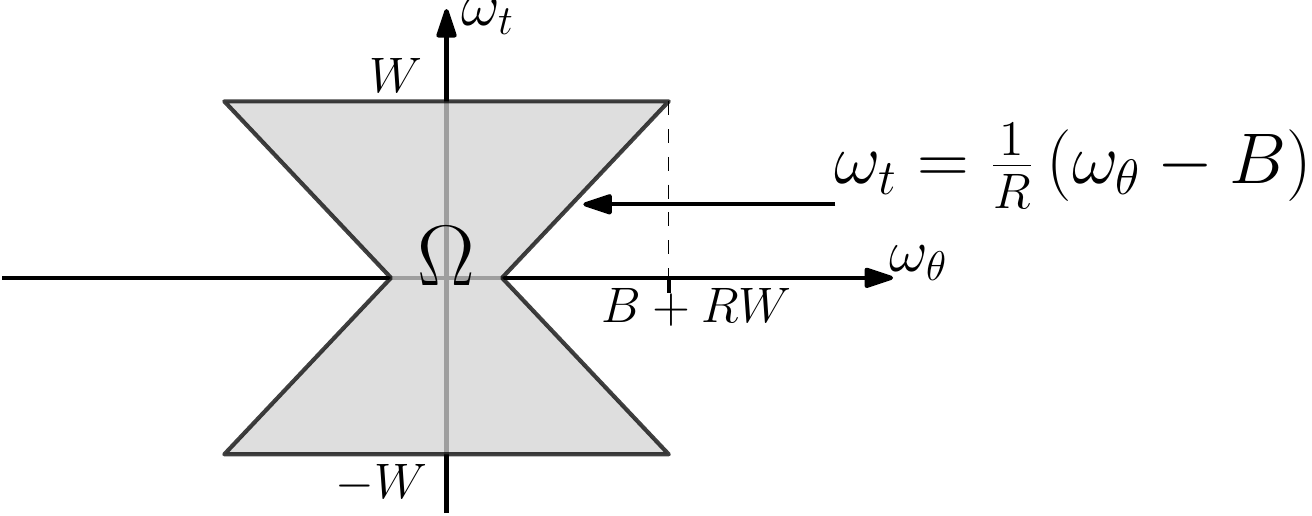}
    \caption{
        The support region of the sinogram in Fourier domain. More than 98\% of the energy is located within this shape when $B=1$. $R$ is the diameter of the scanned object and $W$ is the maximal spatial frequency on the detector axis $t$, determined by the spacing between each adjacent detector $T$.}
    \label{fig:SinogramSupport}
\end{figure}

We now formulate a subspace prior that will aid in the resampling of the PB sinogram to the PP grid, and can also be used for denoising tomographic measurements.
A PB tomographic scan of any bounded object yields a sinogram that has an approximately compact support over its 2D-CSFT $P$~\cref{eq:2D_Sino_CSFT}~\cite{Rattey1981}, in the region
\begin{align}\label{eq:OmegaSS}
   \Omega=\left\{(\omega_t,\omega_\Theta) :
        \left|\omega_{t}\right|  <  W,
        \left|\omega_{\theta}\right|  <  B+\left|\omega_{t}\right|R
        \right\} \,.
\end{align}
Here, $W$ is the maximal spatial frequency of the detector axis, determined by the distance by adjacent detectors $T$, such that $W=1/(2T)$.
The parameter $B$ determines the intersection point with the $\omega_\theta$ axis,
as can be seen in \cref{fig:SinogramSupport}.
The transform of a sinogram $p$ to the 2D-CSFT domain $P$ is given in \cref{eq:2D_Sino_CSFT}.
An illustration of $\Omega$ can be seen in \cref{fig:SinogramSupport}.

The authors in~\cite{Rattey1981} showed that more than $98\%$ of the sinogram's energy is confined to $\Omega$, when the intersection constant is set to $B=1$.
We extend their work by introducing $B$ as a variable parameter in the region $1\leq B$, in order to empirically improve results.
By choosing $1\leq B\leq2$, we get a slightly bigger region $\Omega$, allowing for greater flexibility in our algorithm.
A performance analysis for the selection of different $B$ values, is given in the next subsection.

\begin{figure*}[tb]
    \centering
    \begin{align}
        a(t,\theta) = \begin{cases}
            \frac{2}{\pi}W\left(1+\frac{RW}{2}\right), & t=\theta=0  ,\\
            \big[ 2\theta RW\sin\left(\theta B\right)+\cos\left(B\theta\right) -\cos \left(\theta\left(B+2RW\right)\right)\big]/(2\pi\theta^{2}R), & t=\theta R , \\
            \big[2t\sin\left(Wt\right)\left(B+RW\right)-4R\sin^{2}\left(\frac{Wt} {2}\right)\big]/(\pi t^{2}), & t\neq0,\theta=0 , \\
            \frac{1}{\pi\theta}\left[\frac{\cos\left(Wt-\theta\left(B+RW\right)\right)-\cos\left(B\theta\right)}{t-\theta R}-\frac{\cos\left(Wt+ \theta\left(B+RW\right)\right)-\cos\left(B\theta\right)}{t+\theta R}\right], & \mbox{else}.
        \end{cases}
        \label{eq:TheBigKernel}
    \end{align}
    \rule{\textwidth}{1 pt}
\end{figure*}

The sinogram approximately lies in an SI subspace $\mathscr{A}$, defined by all the functions whose 2D-CSFT is limited to $\Omega$.
Every continuous sinogram $p\in\mathscr{A}$ is spanned by $a(t,\theta)$, so that
\begin{align}
    p \left( t,\theta \right) = \sum_{m,n} \v d[m,n] a\left( t-mT,\theta-n\Theta \right)  \,,
\end{align}
for some (possibly infinite) matrix $\v d[m,n]$.
The kernel $a$ is computed by performing an inverse 2D-CSFT over an indicator function of the region $\Omega$.
\iftoggle{TMI}{
    The analytic expression for $a(t,\theta)$ is shown in \cref{eq:TheBigKernel} and its full derivation is found in the Supplementary Information. 
}{
    The analytic expression for $a(t,\theta)$ is shown in \cref{eq:TheBigKernel} and its full derivation is found in Appendix~\ref{sec:Subspace_Appndx}.
}

Since the region $\Omega$ has finite support in the frequency domain, we deduce that its respective spanning kernel $a$ has an infinite support in the sinogram $(t,\theta)$ domain.
Multiplying the kernel with a window function limits its support, which leads to a low computational burden when performing direct interpolation in the sinogram domain.
In our implementation, we use a Hamming window given by:
\begin{align} \label{eq:Window}
    w \left( t,\theta \right) & =  0.54+0.46\cos\left(2\pi\, \frac{r(t,\theta)}{\sqrt2 K}\right),
\end{align}
with the continuous normalized radius function $r(t,\theta)$ calculated by
\begin{align}
    r \left( t,\theta \right) =  \sqrt{ {\left(\frac{t}{NT} \right)}^2 +
                        {\left(\frac{\theta}{\pi}\right)}^2 } \,. \notag
\end{align}
Here, $K$ is a parameter which sets the effective radius of the window in natural units, and $N$ is the total number of X-ray detectors.

\begin{figure}[t]
    \centering
    \includegraphics[width=\columnwidth]{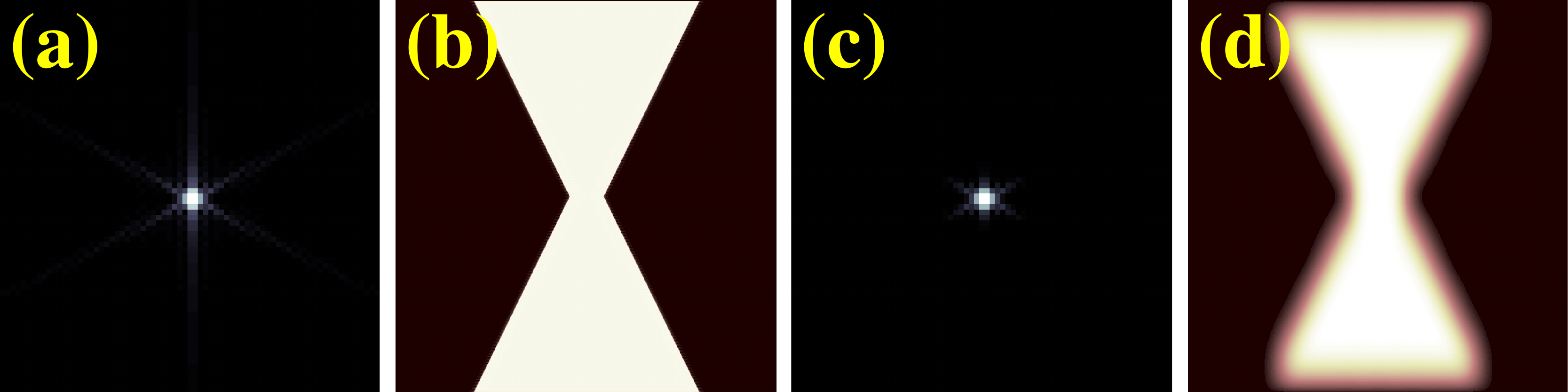}
    \caption{
        (a-b)  The infinite kernel $a(t,\theta)$ that spans $\mathscr{A}$ and its 2D-CSFT.\@
        (c-d) The kernel $a(t,\theta)$ multiplied by the window $w(t,\theta)$ \cref{eq:Window}, with size $K=6$, and its 2D-CSFT.
        This kernel spans the subspace $\mathscr{Q}$ and is denoted by $q$.
    }
    \label{fig:Kernels}
\end{figure}

Once the kernel $a(t,\theta)$ is multiplied by the window $w(t,\theta)$, the sinogram is spanned by a limited number of coefficients at every given $(t,\theta)$ point, such that
\begin{align} \label{eq:p_t_theta_Qb}
    p \left( t,\theta \right) = \sum_{k,l} \v b[k,l] q\left( t-mT,\theta-n\Theta \right)  \,,
\end{align}
where $\v b$ represents the coefficients that span $p$ over the modified subspace $\Q$, with its associated kernel function $q(t,\theta) = a(t,\theta)w(t,\theta)$.
An example of the different kernels $a$ and $q$ (with and without multiplying with a window function, respectively) and their 2D-CSFT is shown in \cref{fig:Kernels} for $K=6$.

\begin{figure}[t]
    \centering
    \includegraphics[width=\columnwidth]{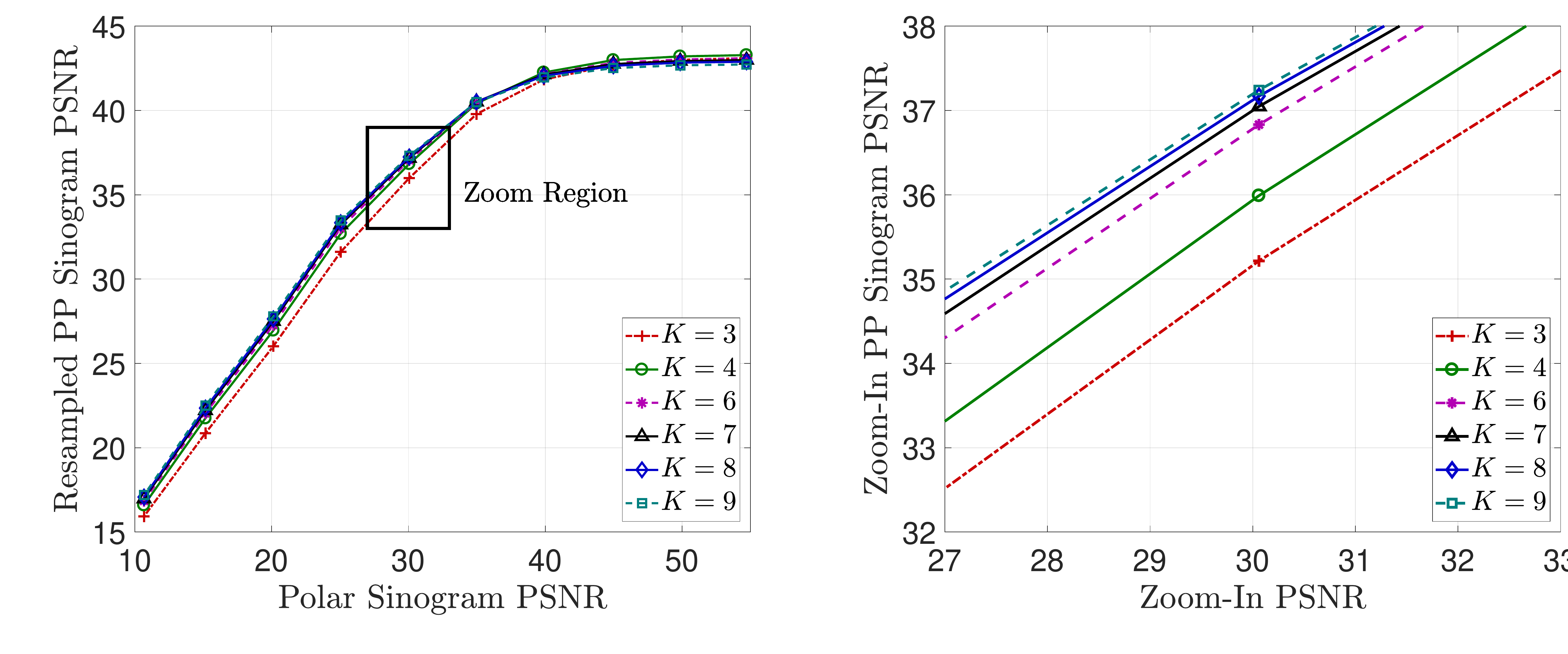}
    \caption{\textbf{Comparison of different window sizes.} We compare the interpolation results for various window sizes $K$ vs.\ the input sinogram's noise ratio. A window size of 6 is sufficiently big for achieving satisfactory results, and the performance of using bigger windows diminishes significantly. The complexity of the resampling procedure is proportional to $K^2$.\label{fig:WindowCompare}}
\end{figure}

The number of coefficients used for a direct computation of a single point using \cref{eq:p_t_theta_Qb} is proportional to $K^2$.
Modifying the choice of the window size $K$ offers the user a trade-off between improving the resampling accuracy and reducing the overall computational complexity.
In our experiments, any value of $K>5$ achieves satisfactory results.
Throughout the rest of the paper we fix $K=6$ to get a good balance between accuracy and speed when implementing the direct convolution computation as described by \cref{eq:p_t_theta_Qb}.
Increasing the size of $K$ further leads to diminishing results, as shown in a simulation study conducted for optimizing the window size parameter, which is presented in \cref{fig:WindowCompare}.


\subsection{Denoising with the Subspace}

We now wish to demonstrate the sinogram denoising capabilities of the derived subspace $\Q$.
The uniformly sampled PB sinogram $\v p$ contains noise induced by various physical phenomena related to the X-ray projection process and properties of the detectors.
By filtering the noisy sinograms with the kernel function $q(t,\theta)$ we cancel the components of $\v p$ that lie outside of the region $\Omega$, therefore removing noise components, while preserving the informative parts of the sinogram required for reconstruction.

To test the subspace denoising performance we first add Poisson noise to a ground truth sinogram $\v p_0$ obtained by simulating a PB projection of a numerical sheep phantom.
The simulation is performed with AIR Tools v1.3~\cite{Hansen2012}.
We denote the contaminated measurements by $\v p_n$.
Next, we denoise the sinograms by convolving with the kernel function $q$ to obtain
\begin{align} \label{eq:p_denoised}
    \v p_\Q = \v p_n * \v q\,.
\end{align}
Here $\v p_\Q$ are the denoised measurements that were filtered by the kernel $\v q$.

The denoised sinogram $\v p_\Q$ is compared to the ground truth sinogram $\v p_0$ for multiple values of $B$, under different noise levels.
The comparison metric is chosen as SNR measured in dB, and calculated by the following function
\begin{align}\label{eq:SNR_Definition}
    \SNR{\v p} &= 20\log_{10} \left\{\frac{\norm{\v p_0}_2}{\norm{\v p - \v p_0}_2}\right\}\,.
\end{align}

\begin{figure}[t]
    \centering
    \includegraphics[width=\columnwidth]{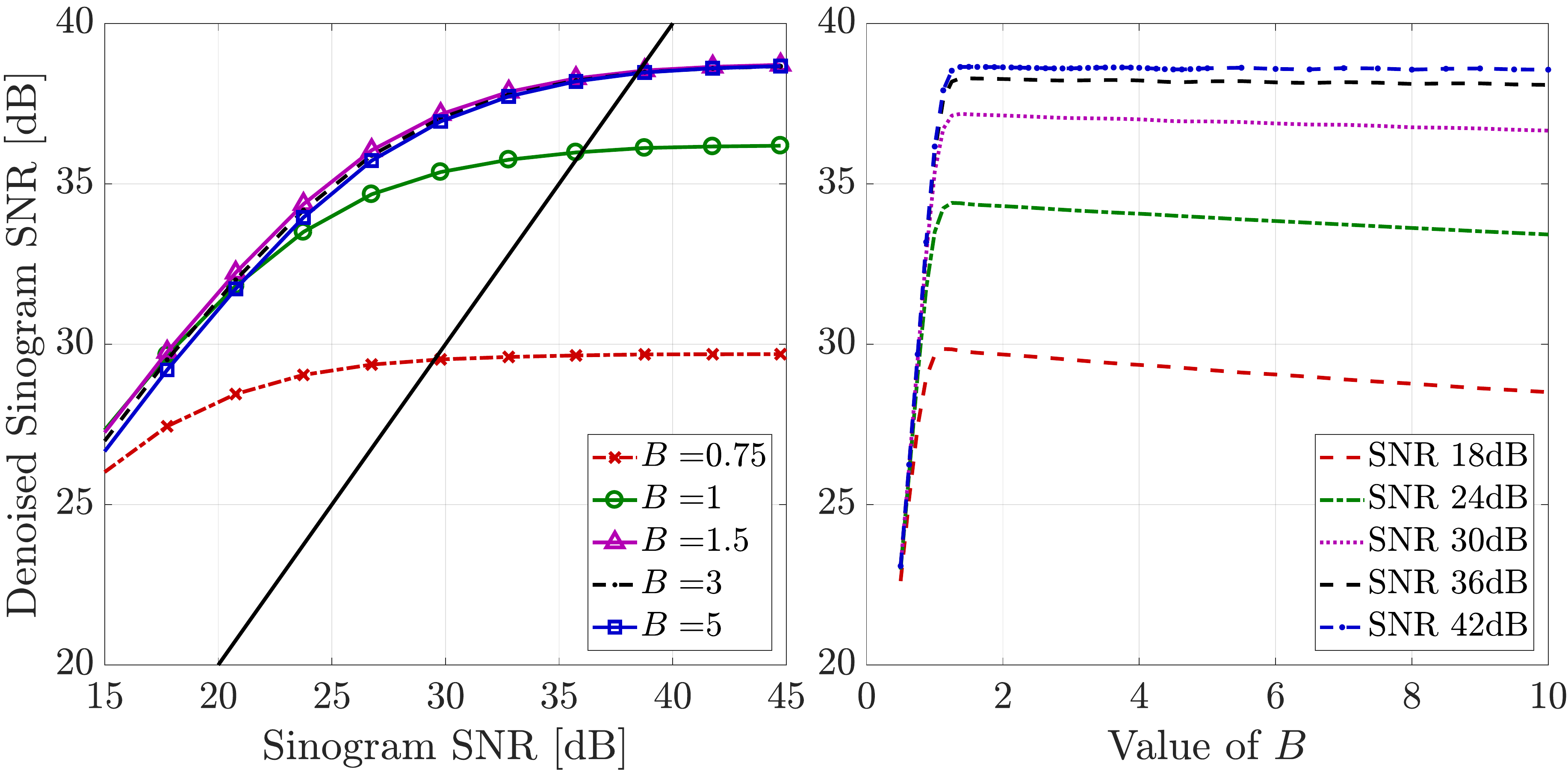}
    \caption{\textbf{Denoising performance and optimal choice for $B$.} 
        In the left graph, the line plots show the output SNR calculated between $\v p_\Q$ and the  $\v p_0$ vs. the input SNR, calculated between $\v p$ and $\v p_0$ for different values of $B$. The diagonal black line is where the input SNR equals to the output SNR.
        In the right graph, the line plots show the output SNR vs. values of B, for different levels of input SNR.\@
        The denoising is effective for moderate levels of SNR, and less effective for low SNR.\@
        In addition, for a relatively clean input there is little to no gain in the denoising procedure, due to the small information loss it incurs.
    }
    \label{fig:B_Denoising_Analysis}
\end{figure}

\Cref{fig:B_Denoising_Analysis} shows that filtering with the subspace kernel successfully denoises the sinogram when compared to the input measured SNR, given by $\SNR{\v p_n}$ .
The black diagonal line denotes equality between the input noise level $\SNR{\v p_n}$ and the output noise level $\SNR{\v p_\Q}$. 
Under moderate noise conditions of $25\db$ SNR, we observe considerable gains of more than $10\db$, when choosing the optimal value for the intersection parameter $B$.
According to the simulation in \cref{fig:B_Denoising_Analysis}, the optimal value of $B$ is empirically found and set to be $B=1.5$.

Note that when trying to denoise images with high levels of SNR (very low noise energy) we observe the performance diminishes.
That is because the region $\Omega$ is an approximation that contains most, yet not all, of the necessary spectral components of the PB radon transform of an object.
In \cref{sec:Reconstruction} we show that in spite of this fact, reconstructing from a denoised sinogram fully maintains fine details of the scanned object.

We conclude that by simply restricting the 2D-CSFT support of a PB sinogram to $\Omega$, we can faithfully denoise it from additive noise.

\subsection{Resampling the Sinogram}

Let $\mathscr{Q}$ be the SI subspace spanned by the windowed kernel $q(t,\theta)=a(t,\theta)w(t,\theta)$.
The discrete PB tomographic samples, $\v p\in\mathscr{Q}$, can then be modeled as
\begin{align}\label{eq:p_Qb}
    \v p       &= \S\Q \v b, \\
    \v p[m,n]  &= \sum_{k,l} \v b[k,l] q\Big( (m-k)T,(n-l)\Theta \Big)\,,\notag
\end{align}
where $\v b$ is the finite dimensional underlying coefficient matrix that represents the continuous sinogram $p(t,\theta)$ over the subspace $\mathscr{Q}$, $\S$ is the PB sampling operator \cref{eq:PB_Sampling} that samples on the uniform PB grid, and $\Q$ is the subspace operator corresponding to convolution with the kernel function $q$.

Our primary goal in this section is to resample the sinogram from PB to PP coordinates.
We first recover the underlying coefficients, denoted by $\v b$, from the PB measurements, $\v p$, by inverting \cref{eq:p_Qb} in the frequency domain.
Once we recover $\hat{\v b}$ we can compute the sinogram over any grid constellation, and specifically over the PP grid.

According to modern sampling theory~\cite{Eldar2015}, if a (continuous) signal $p$ lies within a subspace $\Q$ and sampled by $\S$, we can perfectly reconstruct it by applying the operator $\Q{(\S\Q)}^{-1}$ to its samples $\v p$, as long as the operator $\S\Q$ is invertible.
Otherwise, we can obtain a consistent reconstruction by applying the pseudo-inverse ${(\S\Q)}^\dag$.

In practice, computing the exact inverse ${(\S\Q)}^{-1}$ is a numerically difficult task.
Therefore, we propose to solve the following regularized least squares (LS) convex deconvolution problem,
\begin{align}\label{eq:p-Qb_LS}
    \hat{\v b}=\arg\min_{\v b}\norm{\v p-\S\Q\v b}_{2}^{2}+\rho^2\norm{\v b}^2_2\,.
\end{align}
The operator $\S\Q$ can be written as a linear convolution such that
\begin{align}
    \v p = \v q * \v b\,.
\end{align}
By using the known DFT identity for linear convolution, the above expression is evaluated in $\O(N^2\log N)$ complexity by
\begin{align} \label{eq:DFT_Convlution}
    \v q * \v b = \F_2^{-1}\Big\{\left(\F_2{\v q}\right)\odot \left(\F_2 \v b\right) \Big\}\,,
\end{align}
where $\odot$ denotes an element-wise multiplication between matrices.

The LS objective is regularized with a quadratic norm in order to avoid an ill-posed formulation, which is generally the case when reconstructing from a small number of tomographic measurements.
The $\ell_2$ regularization constant $\rho$ can be arbitrarily small ($10^{-4}$ in our implementation).
By taking $\rho\ra0$ the LS formulation approximates the exact pseudo-inverse ${(\S\Q)}^\dagger \v b$.

The optimal solution to \cref{eq:p-Qb_LS} has a closed form expression given by
\begin{align} \label{eq:b_hat}
    \hat {\v b} & = {\left[(\S\Q)^*\S\Q+\rho^2 \v I \right)}^{-1}(\S\Q)^*\v p,
\end{align}
where $\v I$ is the identity operator.
We note that due to the spatial symmetry of the kernel function $q(t,\theta)$, the operator $\S\Q$ is Hermitian, and the calculation of its adjoint ${(\S\Q)}^*$, is synonymous to applying $\S\Q$. 
We can therefore transform \cref{eq:b_hat} to the frequency domain and calculate it using \cref{eq:DFT_Convlution}.
The LS solution is then computed by element-wise multiplications and divisions on the 2D-DFT of the discrete kernel $\v q$ and PB measurements $\v p$, given by
\begin{align*}
    \hat{\v b}&  = \F_2^{-1} \left\{\frac{(\F_2\v q)\odot(\F_2\v p)}{(\F_2\v q)\odot(\F_2\v q) + \rho^2}\right\}\,.
\end{align*}
Thus, the optimal solution $\hat{\v b}$ to \cref{eq:p-Qb_LS} has a closed form solution, computed in a single step of $\O(N^2\log N)$ complexity.




Once we recover the coefficients $\hat{\v b}$ from the given PB measurements $\v p$, we can resample the continuous sinogram $p$ onto the PP grid, resulting in
\begin{align}\label{eq:InterpolatingToPP}
    \hat{\v p}  & = \Spp\Q{\left(\S \Q \right)}^{\dag}\v p = \Spp\Q\hat{\v b}\,.
\end{align} 
Writing explicitly we have
\begin{align*}
        \hat{\v p}[{m,n}]  &= \sum_{k,l}  \hat{\v b}[k,l]
        q \left( m\hat{T}_n - kT , n\hat{\Theta}_n-l\Theta \right)  \,,
\end{align*}
where $\Spp$ \cref{eq:PP_Spatial_Sampling} is a sampling operator that samples on the PP grid points \cref{eq:PP_Tn_Theta_n}.

We note that the above formulation is completely general and is not limited to point-wise sampling, or a specific subspace.
The sampling operator $\S$ can be generalized so it further takes into account X-ray detector characteristics.

\subsection{Resampling Step Performance}

To assess the performance of our resampling step, we perform a series of simulated tests and comparisons.
A sheep phantom, with resolution of $N=256$, is selected as our scanned object $\v f$.
The digital phantom was computing by scanning a sedated sheep with a Siemens SOMATOM Definition Flash CT scanner.
The projection was performed with the high dosage of 480 mAs over 1250 projection angles, and is therefore chosen to be used as a phantom, due to its relatively high quality.
The PB projection is simulated by the AIR-Tools v1.3~\cite{Hansen2012} package, generating either $60$ or $180$ projections with equispaced angles, each one composed from $256$ detectors.
The synthetic sinogram $\v p$ is then corrupted by Poisson noise with a standard deviation of $\sigma=0.02\ (\max\v{p}-\min\v{p})$ generating the noisy sinogram $\v p_n$.
The interpolation constants are chosen as $B=1.5$, $R=\frac{1}{4\pi}$ and $W=\frac{\pi}{N}$, and the window constant is $K=6$.

\begin{figure}[t]
    \centering
    \includegraphics[width=0.49\columnwidth]{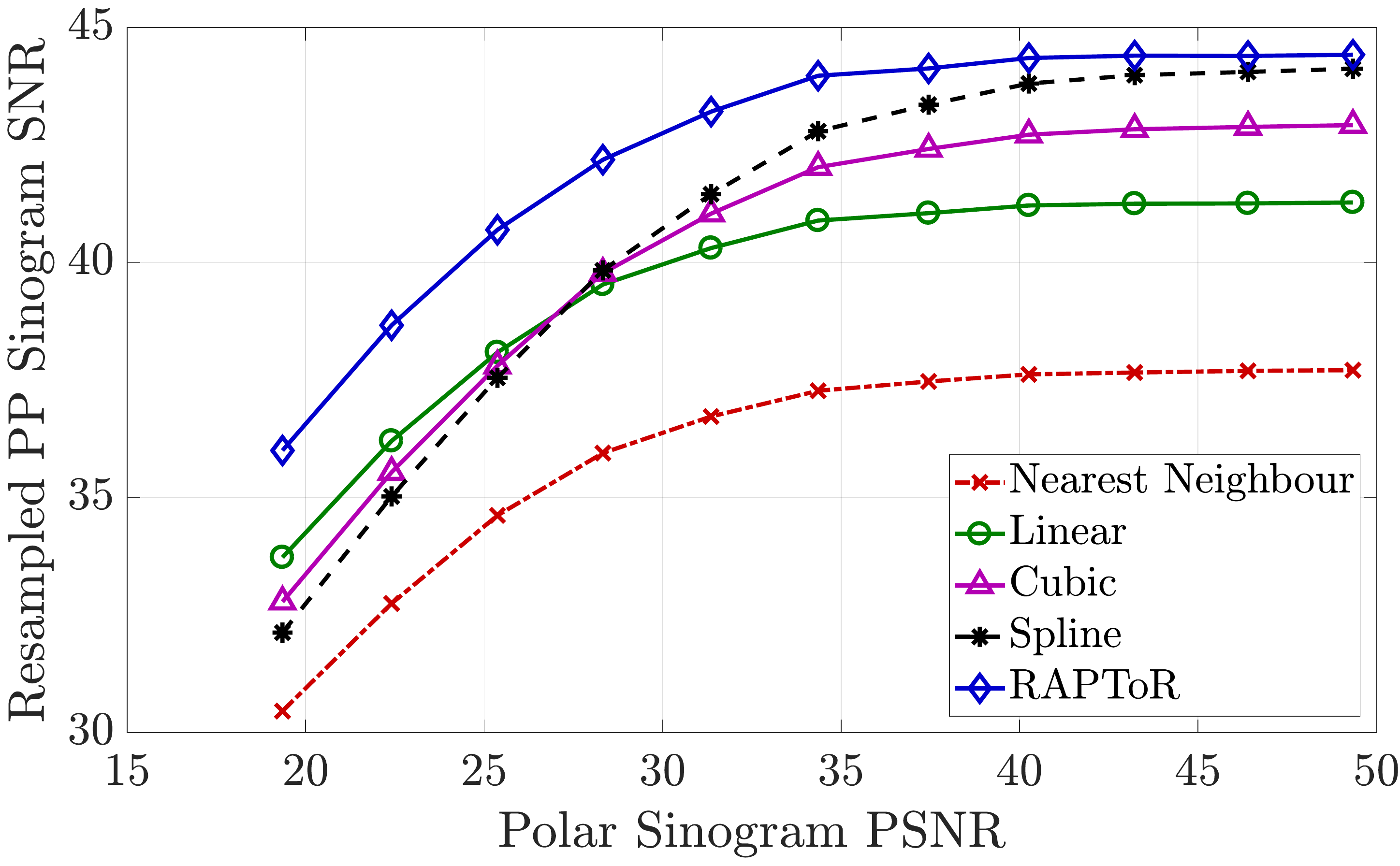}
    \includegraphics[width=0.49\columnwidth]{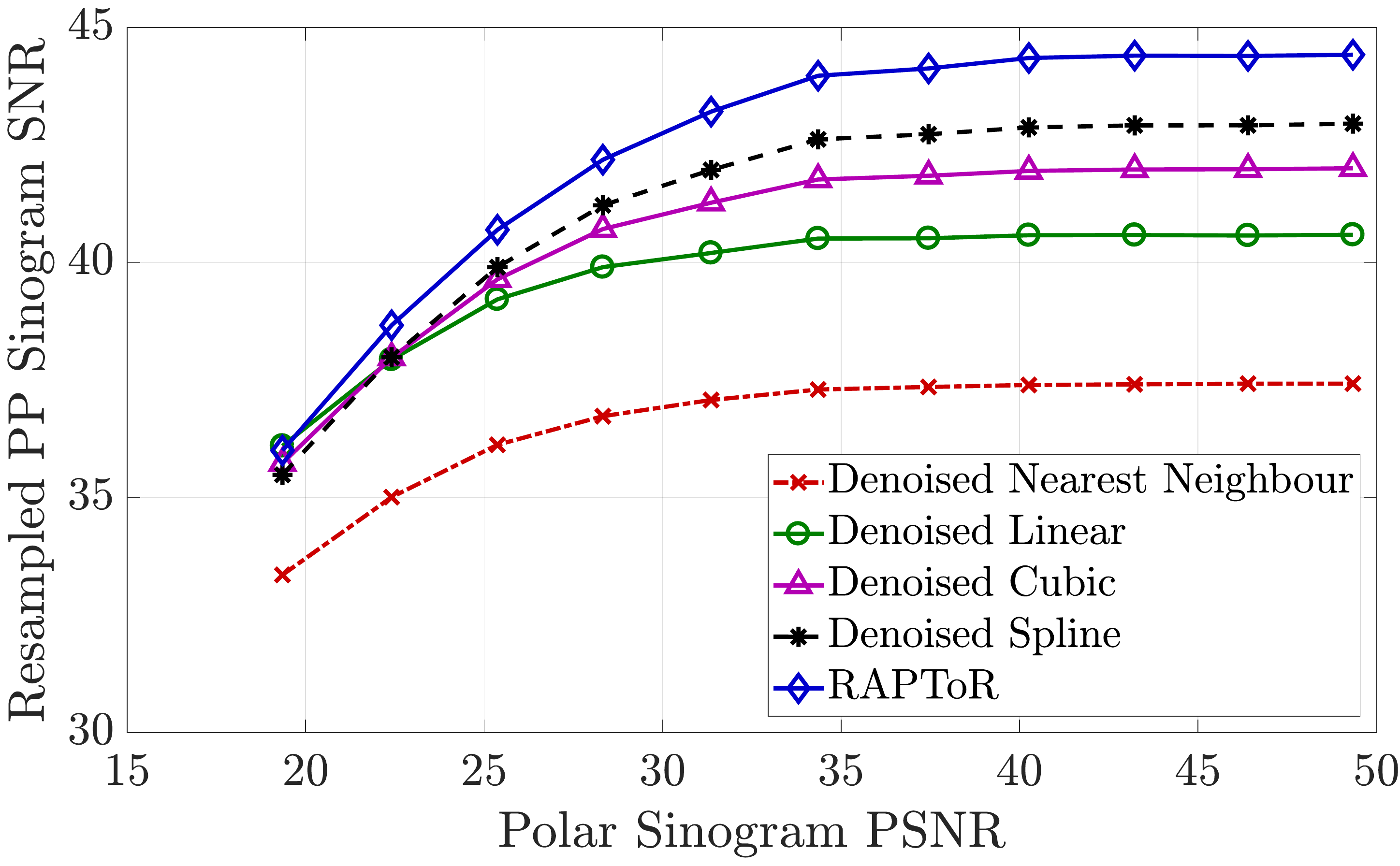}
    \caption{\textbf{Comparison of different interpolation methods.} 
        The horizontal axis depicts the SNR of the inputted polar sinogram.
        The left graph demonstrates the performance of our resampling method compared to other interpolation methods.
        In the left graph we test the performance of our method compared to other approaches, where we use the denoised PB sinogram $\v p_Q$ \cref{eq:p_denoised} as input to the other methods.
        By examining the right graph it is apparent that the performance of our method is not gained only due to the denoising process, but rather by performing the entire resampling algorithm as depicted in \cref{eq:InterpolatingToPP}.
        Our approach achieves superior results for the entire range of input SNR.\@
    }
    \label{fig:Interp_SNR_Compare}
\end{figure}

We compare our method to bilinear, bicubic and a modern spline based interpolation methods, where for reference we use the PPRT \cref{eq:1D_IDFT_p_Sum_s} of the ground-truth phantom.
Each of the regular interpolations are tested for resampling directly from both the noisy PB samples $\v p_n$, and the denoised sinogram $\v p_\Q$ \cref{eq:p_denoised}.
In this way, we can isolate whether our performance is gained by the subspace prior, or due to the generalized resampling algorithm, which solves \eqref{eq:InterpolatingToPP}.
The interpolation results are measured in SNR and structural similarity index (SSIM), displayed in \cref{fig:Interp_SNR_Compare}, confirming that our method successfully denoises the sinogram.
We achieve better results than a two-step approach that first denoises the data by convolving with the kernel \eqref{eq:p_denoised}, and then use a conventional resampling technique.
This confirms that our gains are achieved by both the defined subspace prior, and the generalized resampling approach.


\section{Reconstruction}
\label{sec:Reconstruction}

For reconstructing our scanned object $\hat{\v f}$ from the PP measurements $\hat{\v p}$ we first consider the following TV optimization problem:
\begin{align}\label{eq:PP_Lagrangian}
    \hat{\v f} = \arg\min_f \norm{\hat{\v p}-\Rpp\v f}_2^2+\lambda
    \mathcal{TV}(\v{f})\,,
\end{align}
where $\Rpp$ is the PPRT operator, and $\mathcal{TV}$ is the Total-Variation~\cite{Rudin1992} isotropic norm.
This is a standard optimization approach for reconstruction, shown to produce superior results when coupled with the proposed PP resampling technique, than other state of the art tomographic reconstruction algorithms.

\begin{figure}[t]
    \centering
    \includegraphics[width=0.75\columnwidth]{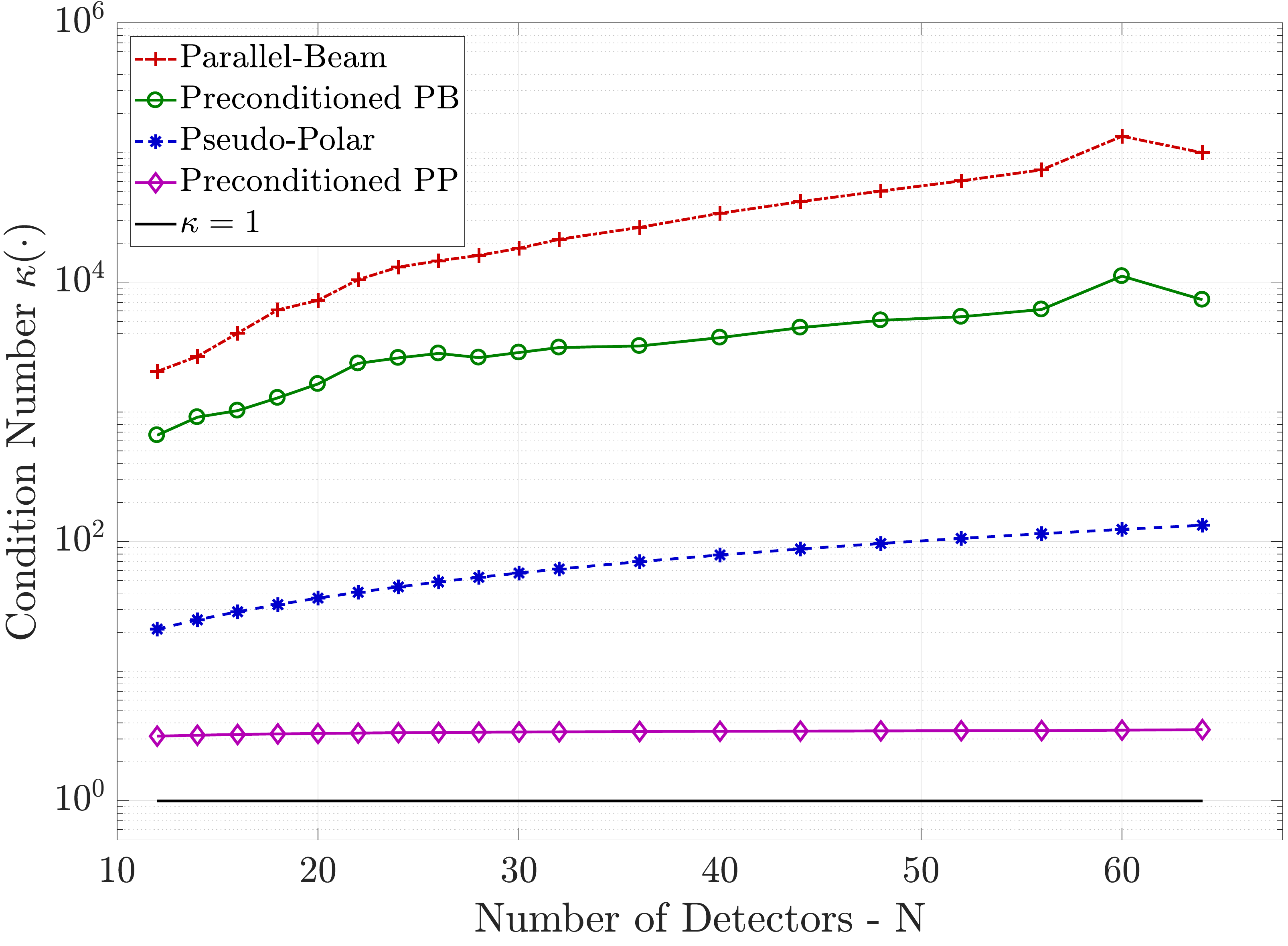}
    \caption{\textbf{Comparison between system condition numbers.}
        The algebraic condition number of various CT measurement systems is plotted as a function of the number of detectors $N$ (logarithmic scale).
        Preconditioning either the PP or the PB system by a diagonal matrix improves the condition number by more than an order of magnitude in each case.
        In either case, the PP system still has a significantly lower condition number.
        A condition number of $1$ is marked by a black line.}
    \label{fig:Condition_Number}
\end{figure}

Before approaching \cref{eq:PP_Lagrangian}, we would like to introduce several modifications to it, that allow for faster and more accurate convergence, with better robustness to noise.
First, the PP operator can be preconditioned by introducing the following operator
\begin{align}
    \M &= \F_1^{-1} M \F_1,
\end{align}
where $\M$ is an element-wise operator described by
\begin{align} \label{eq:M}
    \F_1^*M{\F_1 \v p[m,n]} &= \begin{cases}
        \F_1^{-1}\left[ \frac{1}{6(2N+1)}(\F_1\v p[m,n])  \right], & m=0, \\
        \F_1^{-1}\left[ \left| \frac{m}{2N+1} \right| (\F_1 \v p[m,n]) \right], & \text{else}\,,
    \end{cases}
\end{align}
where $m\in[-N,N]$.
The operator $\M$ is equivalent to applying a high-pass filter to each of the columns in the sinogram it operates on.
\iftoggle{TMI}{
    More information on how we calculate an element-wise preconditioner to the PP system is found in the Supplementary Information.
}{
    More information on how we calculate an element-wise preconditioner to the PP system is found in Appendix~\ref{sec:M_appendix}.
}
When examining the combined operator $\Rpp^*\M\Rpp$, its condition number is significantly lower than that of $\Rpp^*\Rpp$, which leads to faster convergence when using first-order optimization techniques.
A complete comparison between the condition number of the PB and the PP systems, with and without a preconditioner, is shown in \cref{fig:Condition_Number}.
Since the computation of both $\Rpp$ and its adjoint $\Rpp^*$, has complexity on the order of $\O(N^2\log N)$, the overall algorithm complexity remains the same.

Next, we would like to incorporate a statistical model for the measurement noise.
By considering a monochromatic radiation source, the measured X-ray data follows a Poisson distribution~\cite{Elbakri2002a}.
The raw projections measured by the CT detectors are distributed according to
\begin{align}
    \v r[m,n]\sim\text{Poisson}\left\{ I_m e^{-{\v p}[m,n]}+ \epsilon_m\right\} \,,
    \label{eq:Poisson}
\end{align}
where $I_{m}$ is the number of photons emitted towards the $m$th detector, and $\v \epsilon_{m}$ is the average electrical bias in each detector. The line measurements are extracted by calculating the mean of the Poisson distribution in \cref{eq:Poisson}
\begin{align}
    \v p[m,n] = -\ln \left( \frac{\v r[m,n] - \epsilon_m}{I_m} \right) \,.
\end{align}
By assuming that each X-ray path is i.i.d with the others, the authors in~\cite{Elbakri2003} have showed that the first order approximation to the maximum likelihood (ML) estimator that recovers $\v f$ from ${\v p}$ is now given by a weighted LS~\cite{Elbakri2003}
\begin{align}
    \hat{\v f} = \arg\min_{\v f} \norm{ \mathcal{W}\Rpp \v f - \mathcal{W}\hat{\v p} }_2^2
\end{align}
where $\mathcal{W}$ is an element-wise weighting operator dependent on the raw detector measurements $\v r[m,n]$, and given by
\begin{align}
    \mathcal{W}\v p[m,n] = \frac{{(\v r[m,n] - \epsilon_m)}^2}{\v r[m,n]} \v p[m,n].
    \label{eq:W}
\end{align}
Modifying the quadratic norm term in \cref{eq:PP_Lagrangian} together with the operator $\mathcal{W}$, helps to cope better with the additive Poisson noise, and achieve improved reconstructions.

In addition, we would like to reconstruct from a reduced number of measurements. To that end, we define a decimating operator as
\begin{align}
    \mathcal{D}\v p[m,n] = \begin{cases}
        \v p[m,n], & \left\{ m,n \right\}\in\N \\
        0, & \text{else},
        \end{cases}
        \label{eq:D}
\end{align}
where $\mathcal{N}$ is the index set that contains the acquired measurements.

Combining \cref{eq:M}, \cref{eq:W} and \cref{eq:D}, we can define a set of modified operators and a modified system matrix, such that
\begin{align}
    \begin{aligned}
    \A   & = \D\W^{\frac12}\M^{\frac{1}{2}}\Rpp, \\
    \A^* & = \Rpp^*\M^{\frac12}\W^{\frac{1}{2}}\D, \\
    \tilde{\v p}  & = \D\W^{\frac12}\M^{\frac{1}{2}}\hat{\v p},
    \end{aligned}\label{eq:NewOperators}
\end{align}
where $\tilde{\v p}$ contains a limited number of measurements, as described by $\N$.
Since the operators $\D,\ \W$ and $\M$ are all multiplying the sinogram (or its 1D-DFT in the case of $\M$) element-wise, we simply take the square root of each of the multipliers to obtain $\D^{\frac12},\ \W^{\frac12}$ and $\M^{\frac12}$.
Both operators, $\A$ and $\A^*$, exhibit a computational complexity of $\O(N^2\log N)$, since they are composed from a composition of element-wise operators with complexity $\O(N^2)$ and the PPRT, with complexity $\O(N^2\log N)$.

By writing \cref{eq:PP_Lagrangian} with the new operators \cref{eq:NewOperators}, the modified optimization problem is given by
\begin{align}\label{eq:PP_NewLagrangian}
    \hat{\v f} = \arg\min_f \norm{\tilde{\v p}-\A\v f}_2^2+\lambda
    \mathcal{TV}(\v{f})\,.
\end{align}
For solving \cref{eq:PP_NewLagrangian}, we use the FISTA-TV algorithm~\cite{DanielP.Palomar2010,Beck2009b} shown to produce excellent reconstructions from a reduced number of measurements under high noise levels.

The primary gradient step of the solver is given by
\begin{align}
    \v f^{(k+1)} = \v f^{(k)} - \frac{2}{L_f}\mathcal{A}^*\left( \A\v f^{(k)}-\tilde{\v p} \right)\, ,
\end{align}
where $L_f$ is a constant, bigger or equal than the Lipschitz constant of the modified PP system $L$.
This constant is computed by first approximating the largest singular value of $\A$ using the power method.

A summary of the entire algorithm is brought in \cref{alg:MainAlgorithm}, and the FISTA-TV algorithm is presented in \cref{alg:FISTA}.

\begin{algorithm}[t]
    \begin{algorithmic}[1]
        \Statex{\hspace{-1.5 em}\textbf{Input:} A tomographic PB scan, $\v{p} = \mathcal{R} \v f$ \cref{eq:RadonTransform}, the kernel function $\v q[m,n] = (a*w)[mT,n\Omega]$, given in \cref{eq:TheBigKernel,eq:Window}, a TV constant $\lambda$, and a regularization term $\rho$ }
        \Statex{\hspace{-1.5 em}\textbf{Output:} Reconstructed scanned object $\hat{\v f}$}
        \State{\textbf{Find} the underlying coefficients \cref{eq:b_hat}:}
        \Statex{\quad $\hat{\v b}  \gets \F_2^{-1} \left\{(\F_2\v q)\odot(\F_2\v p) / \left[ (\F_2\v q)\odot(\F_2\v q) + \rho^2\right]\right\}  $}
        \State{\textbf{Calculate} the sinogram on the PP grid by}
        \Statex{$\quad\hat{\v{p}} \gets \Spp\Q \hat{\v b}$ \cref{eq:InterpolatingToPP} such that}
        \Statex{$\quad\hat{\v p}[m,n]=\sum_{k,l}  \hat{\v b}[k,l] q \left( m\hat{T}_n - kT , n\hat{\Theta}_n-l\Theta \right)  $ }
        \State{\textbf{Precondition} the input measurements using the operators defined in \cref{eq:W,eq:D,eq:M}:}
        \Statex{$\quad\tilde{\v p} \gets \D\W^\frac12\M^\frac12\hat{\v p}$}
        \State{\textbf{Solve} the following problem using \cref{alg:FISTA}: }
        \Statex{$\quad \hat{\v{f}} \gets \arg\min\nolimits_\v{f} \norm{\v{f-\A\tilde{\v{p}}}}_2^2 +\lambda\mathcal{TV}(\v f) \,,$}
        \Statex{where $\A$ is defined in \cref{eq:NewOperators}}
      \end{algorithmic}
    \caption{PP Tomography Resampling and Reconstruction\label{alg:MainAlgorithm}}
\end{algorithm}

\begin{algorithm}[t]
    \caption{Fast Proximal Gradient Descent for RAPToR\label{alg:FISTA}}
    \begin{algorithmic}[1]
        \Require{} $L\geq L_f$, $\lambda>0$, $K_{\max}$
        \State  Initialize  ${\bf z}_1={\bf x}_0={\bf 0}$, $t_1=1$ and $k=1$
        \For{$k\leq K_{\max}$}
            \State  $\v x_k=\text{prox}_{\mathcal{TV},\frac{2\lambda}{L}}\left(\v x_k-\frac{2}{L}\A^*\left( \A\v z_k - \tilde{\v p} \right)  \right)$
            \State  $t_{k+1}=\frac{1+\sqrt{1+4t_k^2}}   {2}$
            \State  ${\v z}_{k+1}={\v x}_k+\frac{t_k-1}{t_{k+1}}({\v x}_k-{\v x}_{k-1})$
        \EndFor{}
        \State \textbf{return } $\hat{\v f} = {\v x}_{k_{\max}}$
    \end{algorithmic}
\end{algorithm}

\section{Simulation Setup and Results}\label{sec:Setup}

To assess the performance of our algorithm, we perform a series of simulated tests and comparisons.
We divide the tests into two parts, the first compares the quality of the resampling step, and the second evaluates the quality of the scanned object reconstruction.

\begin{figure}[th]
    \centering
    \includegraphics[width=0.49\columnwidth]{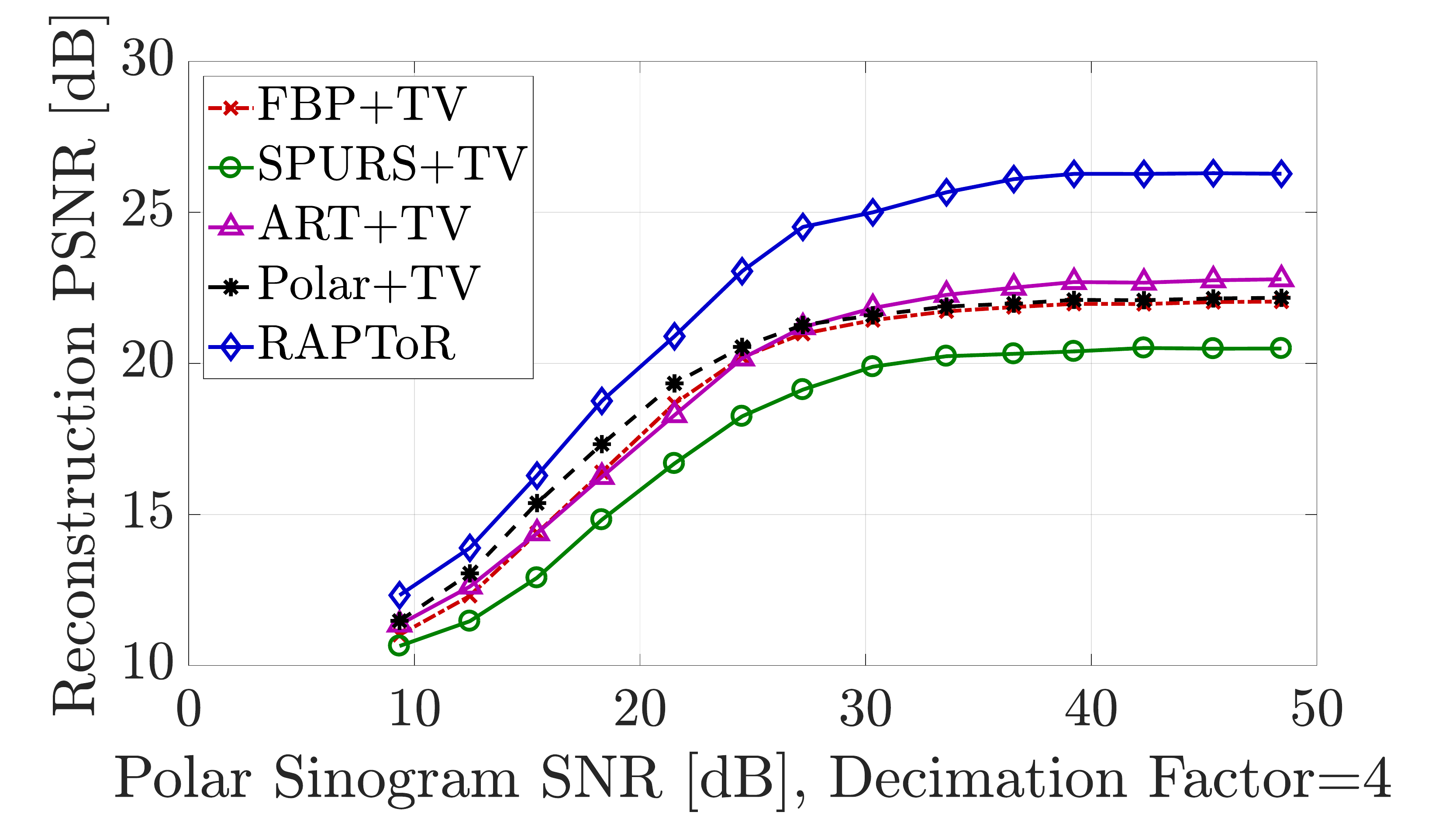}
    \includegraphics[width=0.49\columnwidth]{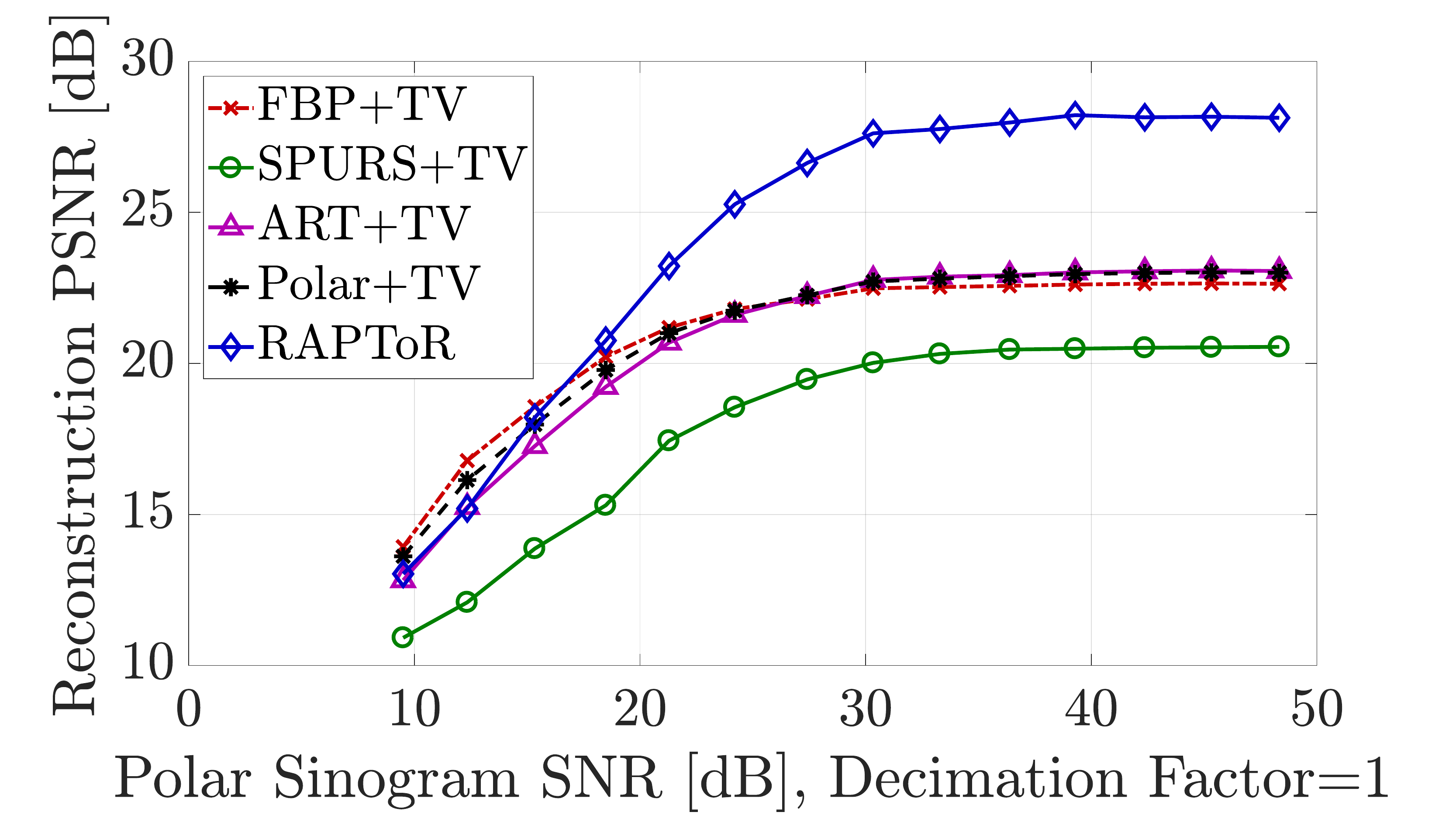}
    \includegraphics[width=0.49\columnwidth]{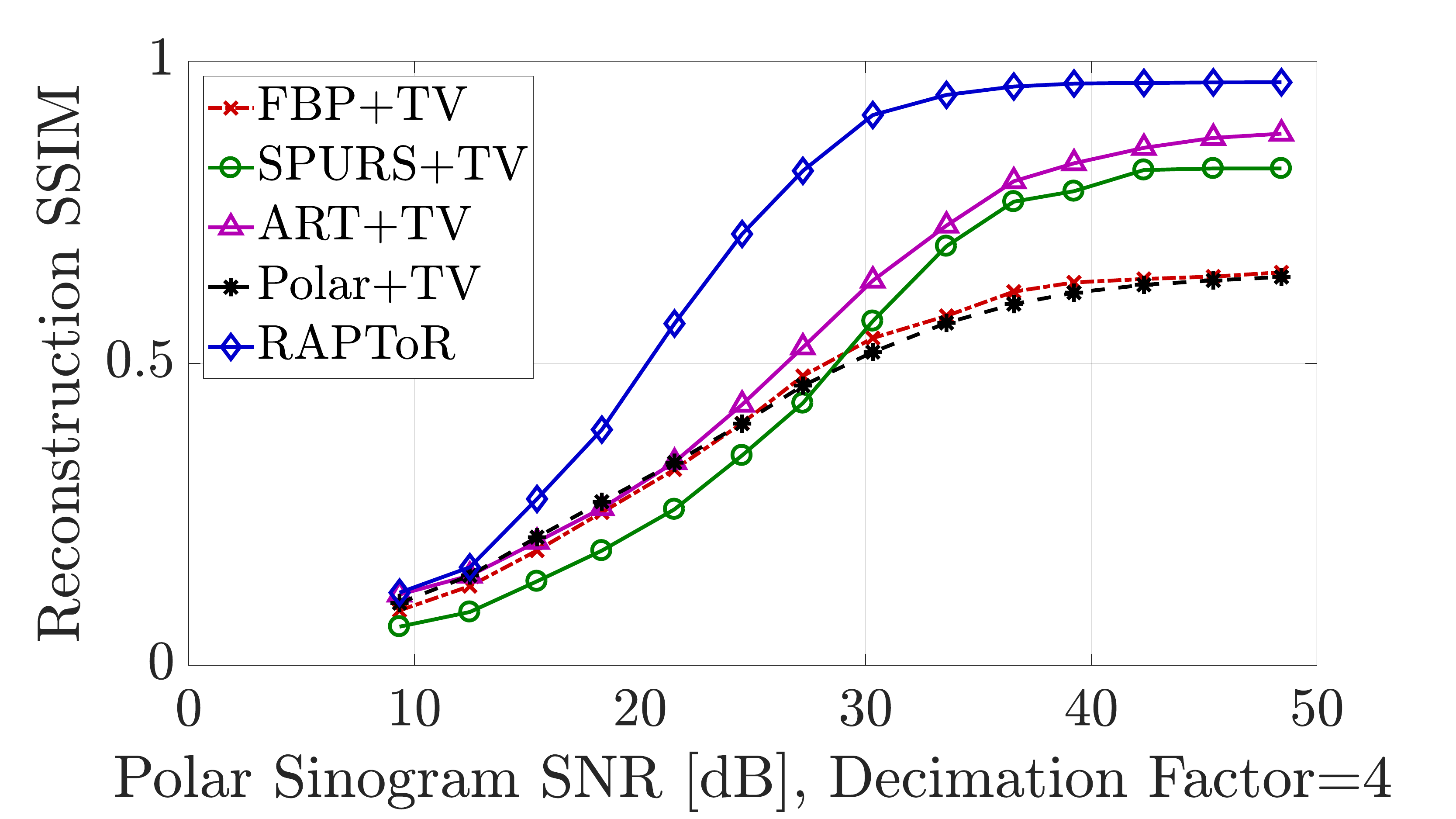}
    \includegraphics[width=0.49\columnwidth]{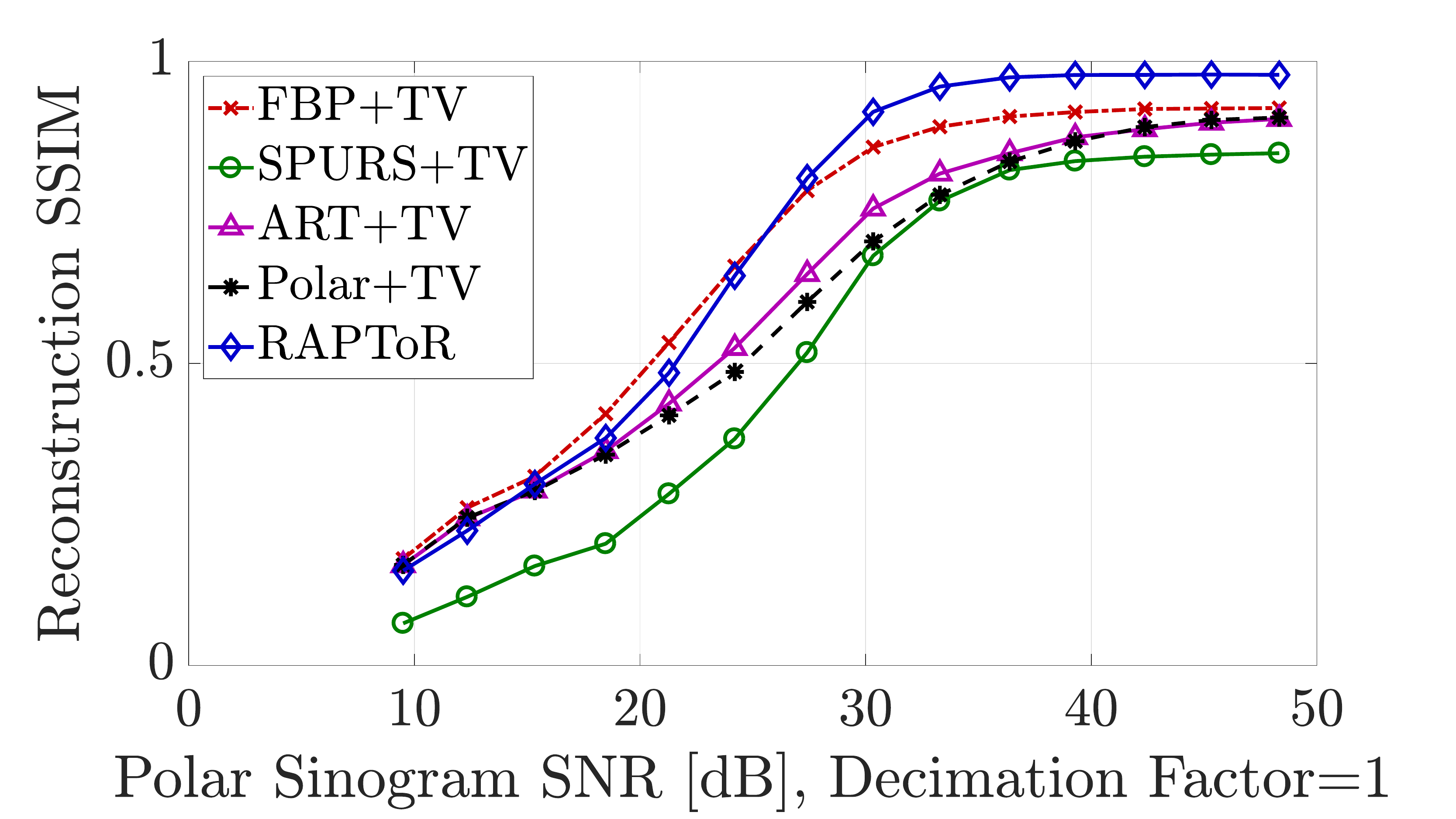}
    \caption{\textbf{Comparison of algorithms for 60 (left) and 180 (right) viewing angles.} The PSNR and SSIM metrics are compared for different state-of-the-art algorithms. Our method, which operates on the transformed PP measurements achieves significantly higher results when the input SNR is greater than $20\db$.} 
    \label{fig:PSNR_SSIM_Recon}
\end{figure}

We select the high dosage sheep phantom (digitized, with resolution $N=256$) as our scanned object $\v f$.
The PB projection is simulated by the AIR-Tools v1.3~\cite{Hansen2012} package, generating only $60$ projections with equispaced angles, each one composed from $256$ detectors.
The synthetic sinogram $\v p$ is then corrupted by Poisson noise with a standard deviation of $\sigma=0.02\ (\max\v{p}-\min\v{p})$.
The interpolation constants are chosen as $B=1.5$ in accordance with \cref{fig:B_Denoising_Analysis}, $R=\frac{1}{4\pi}$ and $W=\frac{\pi}{N}$.

\begin{figure}[t]
    \centering
    \includegraphics[width=0.7\columnwidth]{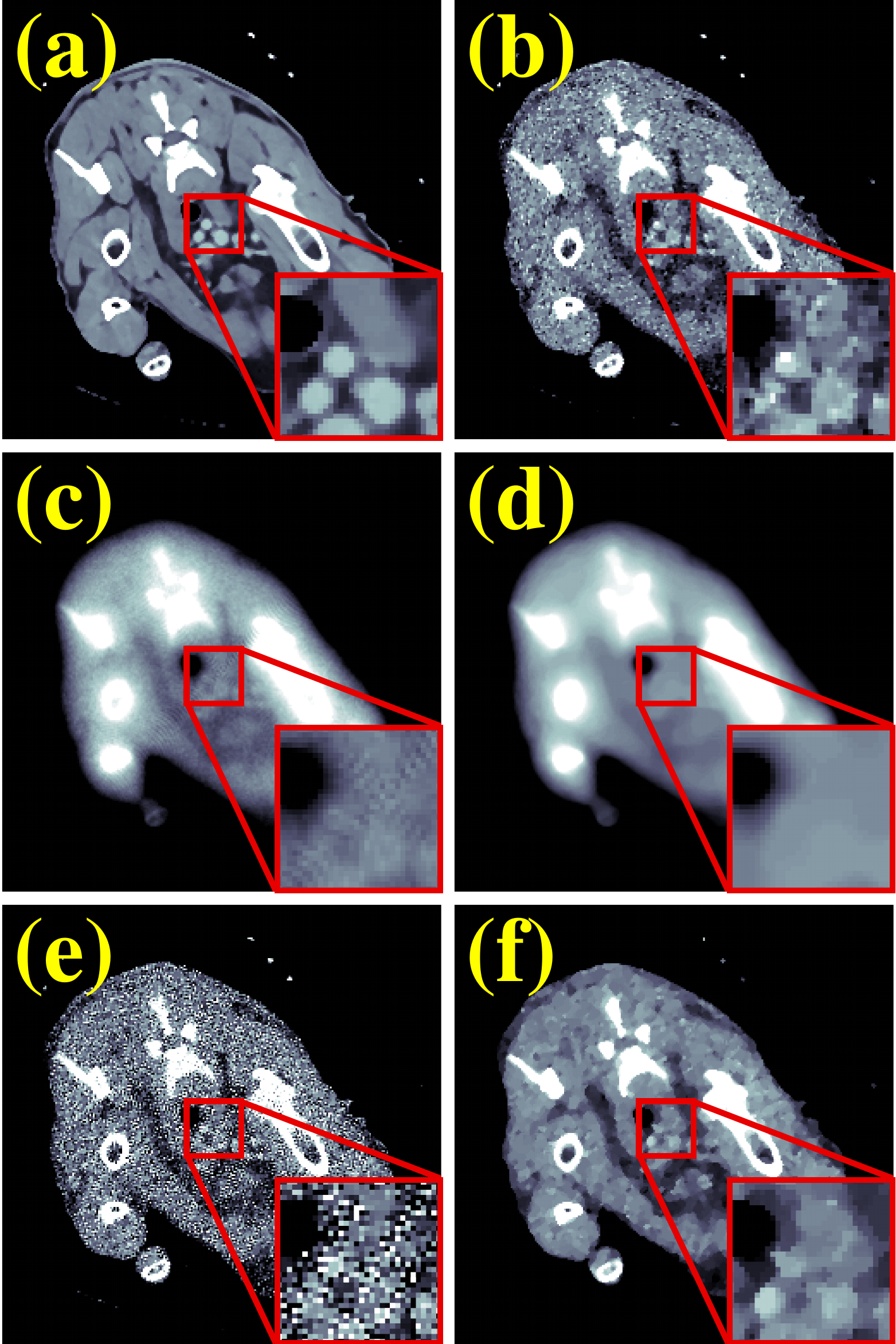}
    \caption{
        \textbf{Comparison of reconstruction results for 180 viewing angles with the same iteration budget.}
        The reconstructions are presented here for the different algorithms.
        The reconstruction is performed from 180 viewing angles, with measurement noise of 28dB.
        (a) ground truth image, (b) FBP followed by TV denoising (FBP+TV), (c) Polar+TV, (d) ART+TV, (e) SPURS+TV, (f) Our method.
        All iterative methods were given the same iteration budget (of 8 iterations), demonstrating similar computational complexity.
    }
    \label{fig:Recons_Same_Iter}
\end{figure}

\begin{figure}[t]
    \centering
    \includegraphics[width=0.7\columnwidth]{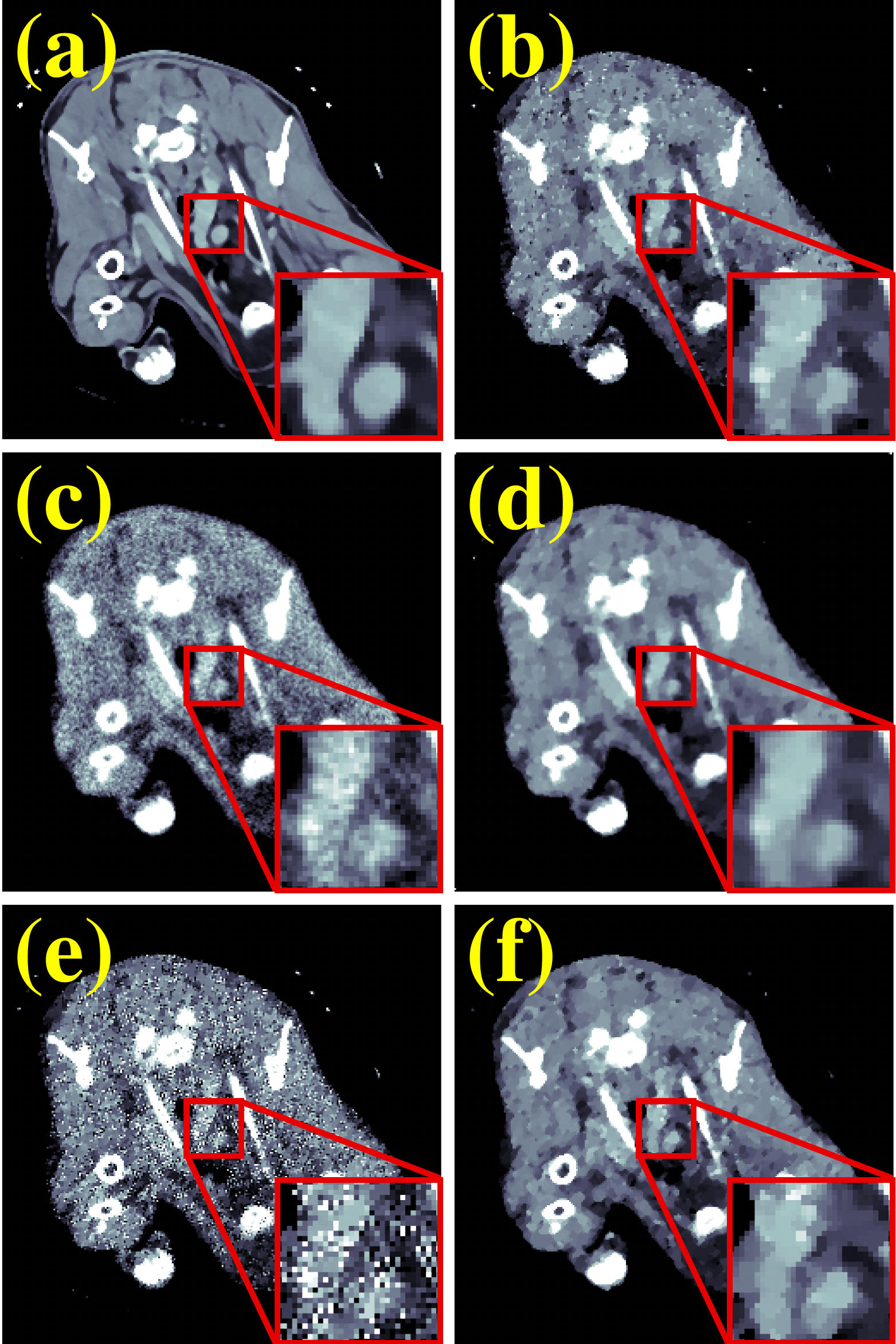}
    \caption{
        \textbf{Comparison of reconstruction results for 120 viewing angles with full convergence.}
        The reconstructions are presented here for the different algorithms.
        The reconstruction is performed from 120 viewing angles, with measurement noise of 24dB.
        (a) ground truth image, (b) FBP followed by TV denoising (FBP+TV), (c) Polar+TV, (d) ART+TV, (e) SPURS+TV, (f) Our method.
        All methods were allowed to run until convergence, specifically Polar+TV and ART+TV used 500 iterations.
    }
    \label{fig:Recons_Conv}
\end{figure}

For the reconstruction procedure we perform a simulated PB projection again with only 60 acquisition angles contaminated with noise, and then resample it to the PP domain using our method.
This scenario is considered challenging for most state-of-the-art solvers available today.
Throughout these simulations the TV coefficient is chosen as $\lambda = 0.03$.
After interpolating to the PP grid, the FISTA-TV algorithm runs for only 20 iterations until converging to the displayed results.
The entire algorithm is compared to FBP, FBP followed by TV denoising, ART+TV~\cite{Chen2008a}, Polar+TV~\cite{Sidky2008b}, and SPURS+TV~\cite{Kiperwas2016}.
The parameters in every algorithm were fine-tuned to run optimally, using an exhaustive grid search for each of them.

For SPURS+TV, we first used SPURS to resample the polar grid samples, given by the 1D-DFT of the PB projections $\F_1\v p$, to the 2D-DFT of the object $\F_2 \v f$, recovering $\hat{\v F}$.
In other words, we resampled in the frequency domain between a polar grid to a Cartesian one.
To make the comparison fair, we then denoised the result by solving the following optimization problem:
\begin{align}
    \v f_\text{SPURS} = \arg\min_{f}\norm{\hat{\v F}-\F_2f}_F^2+\lambda\mathcal{TV}(\v{f})\,,
\end{align}
where an optimal value for $\lambda$ was searched and used.

The final results are shown in \cref{fig:Recons_Conv}.
The other solvers were allowed to run for a significantly higher number of iterations ($>300$), and a longer time until converging to a stable solution.
To demonstrate the relative improvement in time complexity, we conducted an experiment that restricts the iterative methods to the same number of iterations as in our method. The results of this experiment are shown in \cref{fig:Recons_Same_Iter}, where it is evident that competing approaches are still far from a satisfactory solution.
Consequently, we conclude that the time complexity of our combined solution (resampling and reconstruction) is significantly lower than compared solutions.
Due to differences in the hardware implementation of each algorithm, we cannot provide an exact comparison of time complexity, except stating that the asymptotic complexity of our algorithm is slightly smaller than other approaches that has complexity of $\O(N^3\log N)$, while it runs for a number of iterations smaller by more than an order of magnitude.

For testing whether simply denoising the PB sinogram using the subspace prior has a profound effect on the results, we inputted a denoised sinogram $\v p_\Q$ into all other reconstruction approaches. 

\section{Discussion}\label{sec:Discussion}

To conclude, we presented a robust reconstruction scheme for PB tomographic measurements, that performs better than state-of-the-art algorithms with a low computational complexity of $O(N^2\log N)$.
The reconstruction first relies on a resampling process that uses a prior condition on the sinogram, in the form of a subspace, and accurately transforms the PB sinogram to the PP grid, while reducing noise.
The resampling is performed directly on the acquired sinogram measurements, and is shown to reduce noise and achieve high accuracy.
It is shown to perform significantly better than known resampling techniques, either operating in the spatial or frequency domain.
In addition, we showed that a PP-based reconstruction approach, that uses a modern solver which is fed with the resampled PP measurements, produces superior results than state-of-the-art algorithms, with a lower computational complexity and in less time.
The method can also be adapted to modern helical and fan-beam scan methodologies, used in commercial CT scanners, as further discussed in 
\iftoggle{TMI}{the Supplementary Information}{Appendix \ref{sec:Fan_Appendix}}.
We hope this work will pave the way for a CT scanner that uses the PP paradigm for producing clinical images that are taken with significantly less radiation dosage.


\bibliographystyle{Templates/myIEEEtran}
\bibliography{library}

\appendices

\section{Application to Fan-Beam Acquisition}\label{sec:Fan_Appendix}

In this section, we examine the applicability of our method to fan-beam acquistion.
Known algorithms for direct cone-beam or fan-beam reconstruction, namely variants of filtered back-projection (FBP)~\cite{Kak1988c,Hsieh2003} and the Feldkamp-Davis-Kress (FDK)~\cite{Feldkamp1984} algorithms are known to suffer from inaccuracies due to their involved weighting process and interpolations, that compensate for the scan geometry.
These phenomena and others are thoroughly discussed in~\cite{Pan2009,Shaw2014}, where the advantages of iterative reconstruction algorithms are proposed and highlighted.
In practice, due to these limitations, many commercial CTs employ a rebinning step on the measured cone-beam data.
The rebinning process aims to transform the fan-beam projections to an approximately equivalent PB projection, enabling reconstruction by a standard FBP operating on the rebinned measurements.

The problem of accurately rebinning has been already addressed in the past~\cite{Kudo1999,Taguchi2000,Schaller2000,Louis2006}, and scanners today are reconstructing by employing various methods for accurate rebinning.
Many devices use carefully calibrated look-up tables to enable on-the-fly rebinning while the scan is performed to generate PB measurements.
After rebinning the reconstruction process is performed by a FBP.

To further asses the applicability of our method to rebinned PB data, we have conducted a simple experiment.
In the experiment we have generated noisy fan-beam projection of the digital brain phantom~\cite{Alfano2011} phantom, using the AIR-Tools v1.3~\cite{Hansen2012} CT simulator.
We then rebinned the projections to PB equivalent measurements using a simple linear interpolation, and on those measurements we applied our resampling algorithm to transform them to pseudo-polar grid.
For comparison, we reconstructed the object by direct FBP from the fan-beam projections, an FBP used on the rebinned PB projections and a eventually a FBP performed on the pseudo-polar resampled data.

\begin{figure}[t]
	\centering
	\includegraphics[width=\columnwidth]{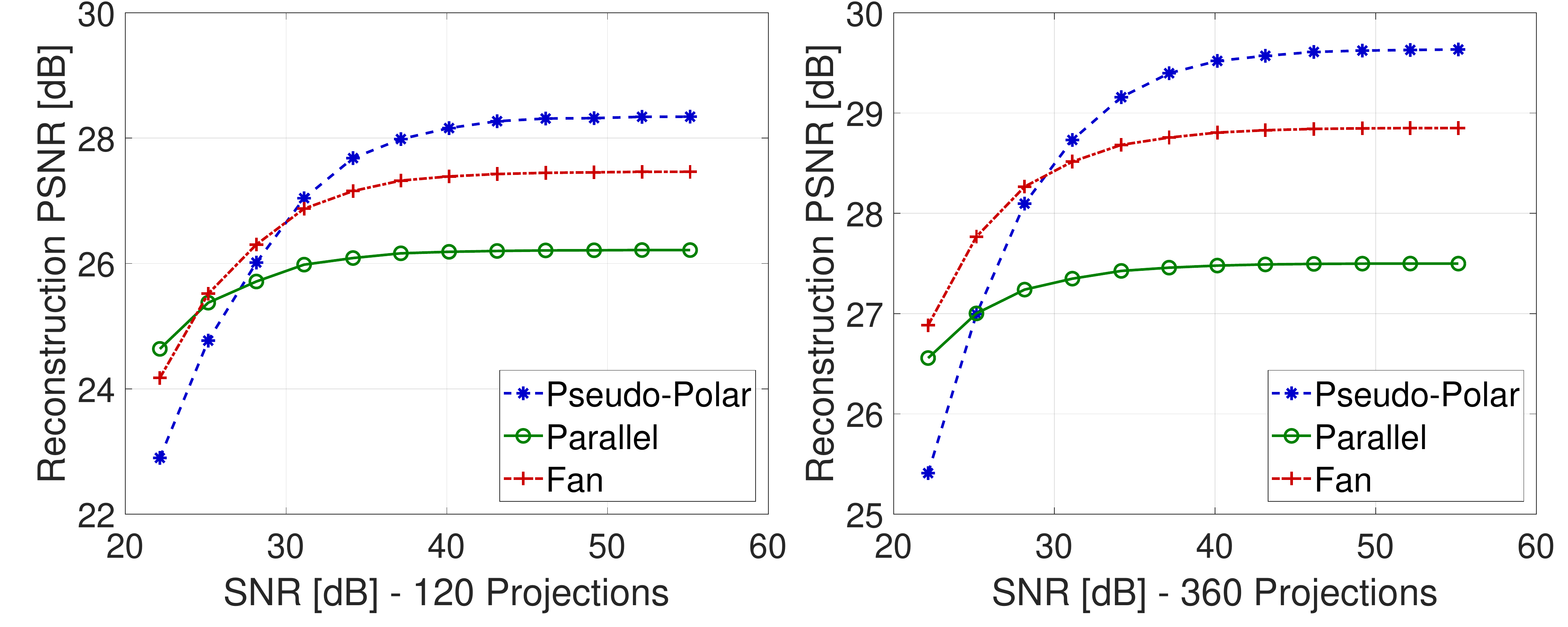}
	\caption{\textbf{RAPToR Performance Comparison.} The numerical phantom is first projected using a fan-beam simulator. Results are compared for direct reconstruction using an FDK~\cite{Feldkamp1984} based algorithm, a parallel-beam FBP performed on rebinned data and an equivalent pseudo-polar FBP performed on data resampled with RaPTOR.
	}
	\label{fig:fanbeam_performance}
\end{figure}

As can be seen in \cref{fig:fanbeam_performance}, our method achieves better reconstructions than both fan-beam and parallel-beam FBP by simply applying the resampling step to the rebinned measurements.
This result is maintained throughout a range of scanning profiles, including different SNR values and the number of angular projections.
After a thorough analysis, we can conclude that even when performing rebinning to PB, which is a lossy transformation, the gains achieved by reconstructing with our method are still superior to other reconstruction approaches. 

\section{Derivation of The Subspace Kernel}\label{sec:Subspace_Appndx}

To find the analytic expression for the kernel $a(t,\theta)$ we perform a 2D inverse continuous Fourier transform over an indicator function of the region $\Omega$ \cref{eq:OmegaSS}:
\begin{align*}
a \left( t,\theta \right)  =& \frac{1}{2\pi}\iint\limits_{\Omega} e^{j \left( t\omega_t+\theta\omega_\theta \right) } \text{d}\omega_t\,\text{d}\omega_\theta \\
=& \frac{1}{2\pi}\dint{-W}{W}{\left(\dint{-B-|\omega_{t}|R}{B+|\omega_{t}|R}{e^{j\left(\theta\omega_{\theta}+t\omega_{t}\right)}}{\omega_{\theta}}\right)}{\omega_{t}} \\
=&
\frac{1}{\pi\theta}\bigg[\frac{\cos\left(Wt-\theta\left(B+RW\right)\right)-\cos\left(B\theta\right)}{t-\theta R}-\\
&-\frac{\cos\left(Wt+\theta\left(B+RW\right)\right)-\cos\left(B\theta\right)}{t+\theta R}\bigg]\,.
\end{align*}

Since this is a continuously differentiable function with respect to $t$ and $\theta$, we can find its limits for each of the irregular points.
For the case where $\theta=0$, the limit is given by
\begin{flalign*}
\lim_{\theta\ra0}a\left(t,\theta\right) &= \lim_{\theta\ra0}\bigg(\frac{2t\sin\left(Wt\right)\sin\left(\theta\left(B+RW\right)\right)}{\pi\theta\left(t^2-\theta^{2}R^{2}\right)}\\
+&\frac{2\theta R\left[\cos\left(Wt\right)\cos\left(\theta\left(B+RW\right)\right)-\cos\left(B\theta\right)\right]}{\pi\theta\left(t^{2}-\theta^{2}R^{2}\right)}\bigg)\\
=&  \frac{2t\sin\left(Wt\right)\left(B+RW\right)-4R\sin^{2}\left(\frac{Wt}{2}\right)}{\pi t^{2}}\,.
\end{flalign*}
The limits for the lines where $t=\pm\theta R$ are given by
\begin{flalign*}
\lim_{t\ra\theta R}&a(t,\theta)= \lim_{\alpha\ra0}a\left(\theta R+\alpha,\theta\right) \\
=&\lim_{\alpha\ra0}\frac{1}{\pi\theta}\bigg[\frac{\cos\left(W\alpha-\theta B\right)-\cos\left(B\theta\right)}{\alpha}- \\
&-\frac{\cos\left(\theta\left(B+2RW\right)\right)-\cos\left(B\theta\right)}{2\theta R}\bigg] \\
=&\lim_{\alpha\ra0}\frac{1}{\pi\theta}\bigg[\frac{-W\sin\left(W\alpha-\theta B\right)}{1} \\
&-\frac{\cos\left(\theta\left(B+2RW\right)\right)-\cos\left(B\theta\right)}{2\theta R}\bigg]\\
=& \frac{W\sin\left(B\theta\right)}{\pi\theta}-\frac{\cos\left(\theta\left(B+2RW\right)\right)-\cos\left(B\theta\right)}{2\pi\theta^{2}R}\,.
\end{flalign*}
For the second transition above we employed L'H\^opital's rule.
The limit for the line $t=0$ is computed by
\begin{align*}
\lim_{\theta\ra0}a\left(t,\theta\right)&=\lim_{\theta\ra0}\frac{-2}{\pi\alpha}\left[\frac{\cos\left(\theta\left(B+RW\right)\right)-\cos\left(B\theta\right)}{\theta R}\right]\\
&=\lim_{\theta\ra0}\frac{1}{\pi\theta}\left[\frac{\theta^{2}{\left(B+RW\right)}^{2}-B^{2}\theta^{2}}{\theta R}\right]\\
&=\frac{2}{\pi}W\left(1+\frac{RW}{2}\right)\,.
\end{align*}
For computing the limit of the point $(t,\theta)=(0,0)$ we choose the trajectory $(t,\theta)=(0,\alpha)$ when $a\alpha\ra 0$:
\begin{align*}
\lim_{t,\theta\ra0}a\left(t,\theta\right)=&\lim_{\alpha\ra0}\frac{-2}{\pi\alpha}\left[\frac{\cos\left(\theta\left(B+RW\right)\right)-\cos\left(B\theta\right)}{\theta R}\right]\\=&\lim_{\theta\ra0}\frac{1}{\pi\theta}\left[\frac{\theta^{2}{\left(B+RW\right)}^{2}-B^{2}\theta^{2}}{\theta R}\right]\\
=&\frac{2}{\pi}W\left(1+\frac{RW}{2}\right)\,,
\end{align*}
thus concluding the derivation of the subspace kernel. Its derived expression is given in \cref{eq:TheBigKernel}.

\section{The Pseudo-Polar Preconditioner}
\label{sec:M_appendix}

We would like to find an algebraic preconditioner to the PP system, that can be computed with low computational complexity.
The ideal algebraic preconditioner to the PP system is given by the following expression
\begin{align}
\hat \M = {\left(\Rpp\Rpp^*\right)}^{-1}\,.
\label{eq:M_Ideal}\tag{A1}
\end{align}
This preconditioner (if exists), leads to a condition number of $1$ for the pseudo-polar tomographic system, since by plugging it in we get:
\begin{align*}
\Rpp^* \hat \M \Rpp = \Rpp^*\left(\Rpp\Rpp^*\right)^{-1}\Rpp=\v I, 
\end{align*}
where $\v I$ is the identity operator.

\begin{figure}[t]
	\centering
	\includegraphics[width=\columnwidth]{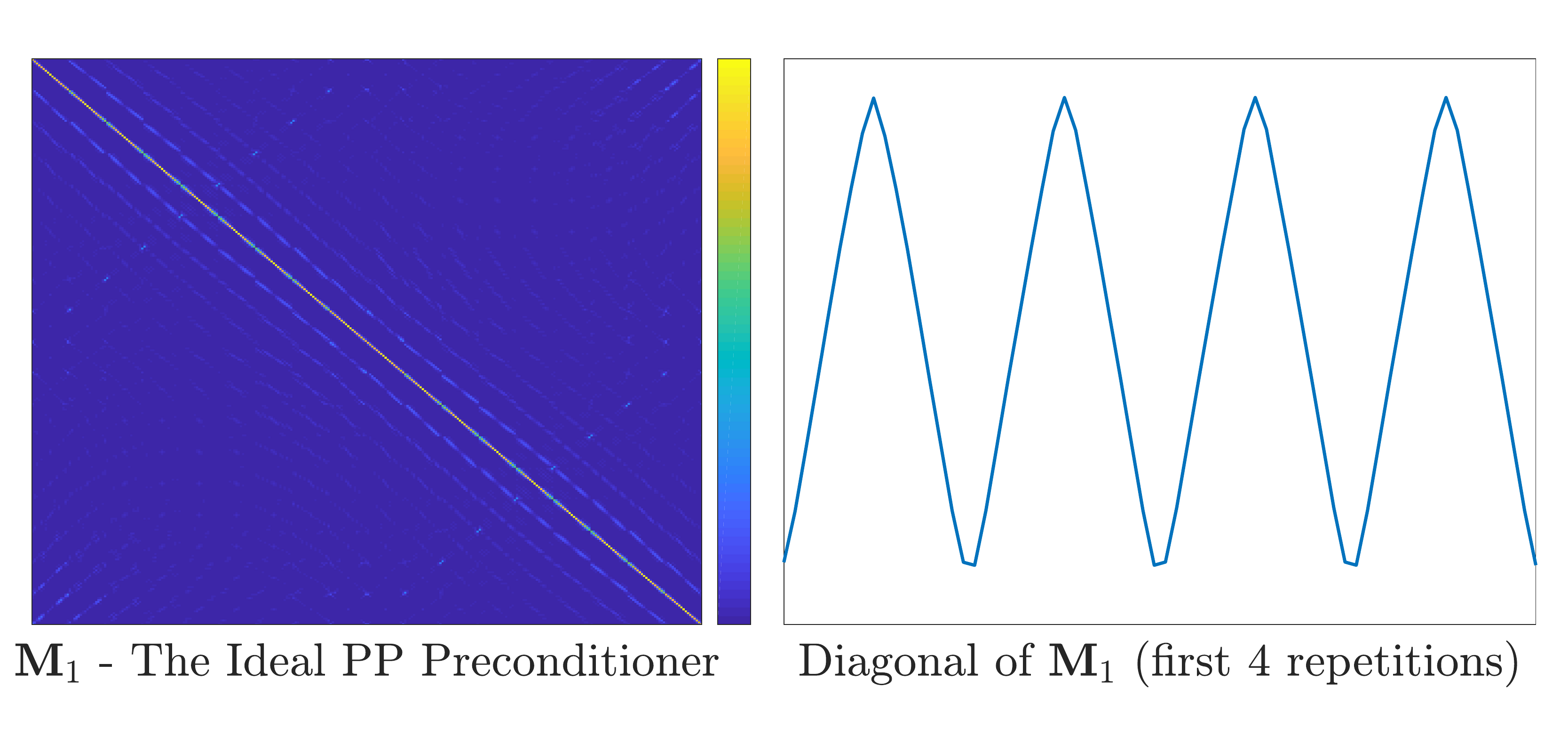}
	\caption{
		The optimally calculated preconditioner (after applying the operator $\F_1$ on it) is shown on the left for a small dimensional PP transform ($N=16$). It was calculated by numerically evaluating \cref{eq:M_Ideal}. More than $\%90$ of the energy is found on the main diagonal. On the right, for clarity we plotted only the first 4 repetitions from the diagonal of $\hat \M$, showing that the preconditioning is the same for every column in the PP sinogram. For bigger dimensions of $N$, we apply an element-wise version of the operator $\M$, as depicted in \cref{eq:M}.}
	\label{fig:M_Diagonal}
\end{figure}

For finding a preconditioner that efficiently approximates $\hat\M$, we first numerically evaluate the following matrix inversion problem:
\begin{align}
\label{eq:M_0_Definition} \tag{A2}
\v M_0 = \left( \v R_{\text{pp}}\v R_{\text{pp}}^*+\varepsilon\v I\right)^{-1}\,.
\end{align}
Here $\v R_{\text{pp}}\in\R^{(2N+1)^2\times N^2}$ and $\v R_{\text{pp}}^*\in\R^{N^2\times{(2N+1)}^2}$ are the matrix form of the PPRT operators $\Rpp$ and its adjoint, $\v I$ is the identity matrix, and $\varepsilon$ is a regularization constant chosen arbitrarily small ($1\text{e}^{-5}$ in our case), that makes sure the inverse to \cref{eq:M_0_Definition} uniquely exists.
These operate on vectorized scanned objects $\text{vec}\{\v f\}\in\R_{+}^{N^2}$ and PP sinograms $\text{vec}\{\v p\}\in\R^{(2N+1)^2}_{+}$ accordingly.
In practice, due to large dimensions of the matrix versions of the PP operators, the ideal preconditioner $\v M_0$ can only be evaluated for small dimensions.
In our experiment, it was computed for object dimensions of $32\times32$ pixels ($N=32$).

Even if we could compute and store the ideal preconditioner for conventional dimensions, applying it in every iteration is a prohibitive step with computational complexity of $\O(N^4)$.
In an attempt to alleviate this computational burden, we apply the 1D-DFT operator to $\v M_0$, such that
\begin{align*}
\v M_1 = \F_1 \v M_0.
\end{align*}

Most of the energy in $\v M_1$ is concentrated on its main diagonal, as can be seen in \cref{fig:M_Diagonal}.
In addition, by plotting the main diagonal of $\M_1$, we observe it is periodic with a periodicity equal to the number of X-Ray detectors of the examined system.
We can therefore approximately substitute the ideal preconditioner $\v M_0$ by taking a diagonal matrix, with its main diagonal taken from $M_1$ and applying an inverse 1D-DFT operator $\F_1^{-1}$ on it.
This procedure corresponds to filtering each of the PP projections for every acquisition angle with the same filter, and has computational complexity of $\O(N^2\log N)$.

To summarize, we numerically show that in the case of the PP system, a good approximation to the optimal algebraic preconditioner is given by a 1D filter that operates on the columns (detector axis wise) of the sinogram.
The exact filter can be deduced by examining and extracting the diagonal of $\v M_1$ for small dimensions, for which its exact computation is numerically feasible.



\end{document}